\documentclass[twoside]{IEEEtran}
\usepackage{cite}
\usepackage{amsmath,amssymb,amsfonts}
\usepackage{algorithm,algorithmic}
\usepackage{graphicx}
\usepackage{textcomp}
\usepackage{psfrag}
\usepackage{array}
\usepackage{subfigure}
\usepackage{color,xspace}
\usepackage{myabbrev} 
\usepackage{fancyhdr}

\newcommand{\CC}{\mathcal{C}}

\begin{document}

\title{Asynchronous Distributed Event-Triggered Coordination for Multi-Agent Coverage Control}
\author{Mohanad~Ajina,
        Daniel~Tabatabai,
        and~Cameron~Nowzari
\thanks{Manuscript submitted April 24, 2019; submitted October 07, 2019.}
\thanks{M. Ajina, D. Tabatabai and  C. Nowzari are  with the department of Electrical and Computer Engineering, George Mason University, 4400 University Drive, Fairfax, Virginia 22030. email: \texttt{\{majina,dtabatab,cnowzari\}@gmu.edu}}}

\markboth{IEEE TRANSACTIONS ON CYBERNETICS}{Ajina \MakeLowercase{\textit{et al.}}: Asynchronous Distributed Event-Triggered Coordination for Multi-Agent Coverage Control}

\maketitle

\begin{abstract}
This paper re-visits a multi-agent deployment problem where agents are restricted from requesting information from other agents as well as sending acknowledgments when information is received. These communication constraints relax the assumptions of instantaneous communication and synchronous actions by agents (request and response actions).
In this paper, we propose a fully asynchronous communication aware solution to the multi-agent deployment problem that uses an event-triggered broadcasting strategy. Unlike all existing triggered solutions, our event-triggered broadcasting algorithm relies on agents to decide when to broadcast (push) information to others in the network without the need for a response from other agents. In addition, the proposed strategy determines how best to move when up-to-date information is unavailable and cannot be requested. The algorithm is capable of achieving similar levels of performance to that of a continuous or periodic strategy. Our solution is proven to achieve asymptotic convergence and simulation results are provided to demonstrate that the proposed event-triggered broadcasting algorithm can achieve an adequate level of performance under the communication constraints.
\end{abstract}

\IEEEpeerreviewmaketitle

\section{Introduction}
\IEEEPARstart{R}{esearchers} have gained increasing interest in Wireless Sensor Networks (WSNs) due to the variety of applications that can benefit from their exploitation. A WSN is a collection of sensors that have the ability to communicate with one another through a shared wireless spectrum (SWS), and the use of WSNs have been observed in a variety of indoor and outdoor monitoring, tracking and security applications~\cite{AP-JS-DW:04,RZ-DY-YW:07,IFA-WS-YS-EC:02,SHY:14}. They can be constructed with a small number of sensors, such as in home security applications~\cite{NN-TH-VR-GV-EM-GM-AM-TP:00}, or can be constructed with several thousand sensors as in environmental monitoring applications~\cite{JKH-KM:06}. Sensors can be equipped with locomotion capabilities allowing them to reconfigure based on changes to the environment as well as allowing them to be deployed to areas that would otherwise be deemed infeasible for placement. This paper is interested in optimally deploying mobile sensors with the use of a practical communication model that does not require {instantaneous communication and synchronous actions by mobile sensors} while performing the deployment task.

As the number of sensors in a WSN increases, the amount of communication between sensors increases as well. This rise in communication comes with a cost where more messages that are exchanged between sensors can create a bottleneck over the SWS. This can increase the frequency of packet drops and can create transmission delays in the network. Current data suggests that the SWS is already overcrowded and by 2030 the demands on SWS applications will be 250 times greater than present day~\cite{SC:16}. 
{Previous efforts in~\cite{CN-JC:11-auto,CN-JC:15-crc,ZZ-XZ-LJ:18,NH-YM-ST:15} have been made to reduce the amount of communication between sensors while still maintaining an adequate level of performance. This is achieved using event-, self- or team-trigger strategy to communicate with the requirement of instantaneous communication and synchronous actions by mobile sensors. In these works, when a trigger has occurred, an mobile sensor must request information, and the neighboring mobile sensors must respond to that request immediately.
We investigate this problem by relaxing the assumption that mobile sensors must communicate instantaneously and take synchronous actions. In other words, the existing works rely on the fact that the mobile sensor can obtain information when it deems necessary. In contrast, our work builds on mobile sensors pushing (broadcasting) information to others, which means that mobile sensors can no longer request information when they desire.
Our solution is a fully asynchronous event-based broadcasting communication strategy that further reduces the number of messages exchanged between mobile sensors.
For the remainder of this paper, we more generally refer to mobile sensors as agents.}

\paragraph*{Literature review}
Similar to previous works that study the multi-agent deployment problem, we make use of \textit{Voronoi partitions} in order for the agents to use distributed gradient descent laws to converge to the set of locally optimal solutions~\cite{JC-SM-TK-FB:02-tra}. Also, the authors in~\cite{AMCS-YY:05} utilized Voronoi partitioning techniques to verify whether or not a set of sensors sufficiently covers the environment. The authors assume that the agents always have access to all other agents' locations at all times. In~\cite{ AB-XF:07, GH-ST-YT-HN:12}, the assumption of always having access to other agent positions is relaxed whereby deployment to the set of critical points is achieved based only on local information provided by an agents' neighbors.  
However, all of the above-mentioned works assume that continuous or periodic communication occurs between agents. It is the desire of this paper to relax the requirement of continuous or periodic communication {without the need for instantaneous communication and synchronous actions by agents}. To do so, we turn to the active research areas involving distributed event-, self-~\cite{WPMHH-KHJ-PT:12} and team-triggered strategies~\cite{CN-JC-GJP:15-acc}.

The principal idea behind event-, self-, and team-triggered strategies is that they provide a means for which agents are capable of only communicating when a designed triggering criterion occurs. The usual consequence of these methods is that agents communicate aperiodically and therefore less frequently than continuous or periodic strategies~\cite{YF-GF-YW-CS:13, CN-JC-GJP:15-acc,CN-JC:16-auto,XP-ZL-ZC:16,GSS-DVD-KHJ:13,MA-CN:18,CN-JC:11-auto}. 
Consider an agent that continuously or periodically communicates its state to neighboring agents. This may be considered a waste of resources. Instead, a more efficient approach would be for the agent to only communicate its state information when it anticipates something may have changed or when something may have gone wrong. 

A brief explanation of the differences between the three mentioned triggering strategies follows. An event-triggered strategy requires an agent to monitor its state until some conditions are met to initiate communication between agents. A self-triggered strategy requires an agent to determine when the next triggering time will occur given its current information. Finally, a team-triggered strategy requires a group of agents to collectively determine when to initiate communication by using an event- and/or self-triggered strategy \cite{CN-JC:15-crc}.
{Generally speaking, we would like to note that there is no clear advantage of one strategy over the others, and depending on the problem or application, one strategy maybe suitable where the others are not.}

All triggering strategies in~\cite{YF-GF-YW-CS:13, CN-JC-GJP:15-acc,CN-JC:16-auto,XP-ZL-ZC:16,GSS-DVD-KHJ:13,MA-CN:18,CN-JC:11-auto} have been shown to reduce the amount of message exchange and/or the amount of communication power consumed between the agents in a number of different multi-agent coordination tasks. This is due to the fact that agents only communicate when there is a need.
In~\cite{YF-GF-YW-CS:13}, an event-triggered algorithm for a multi-agent rendezvous task was developed that reduces the communication between the agents. For multi-agent average consensus tasks, self-triggered algorithms in~\cite{CN-JC:16-auto,XP-ZL-ZC:16} and an event-triggered broadcasting algorithm in \cite{GSS-DVD-KHJ:13} have been shown to reduce the amount of communication {and the amount of communication power consumed during the task.} 
In~\cite{XL-JS-LD-JC:17,MZ-CP-WH-YS:17}, the communication power consumed was shown to be reduced for a leader-follower consensus problem using an event-triggered algorithm.
The most relevant to this paper is the work involving event-, self- or team-triggered strategies applied to the deployment problem. 

In relation to multi-agent deployment, triggering strategies have been shown to reduce the amount of communication as in~\cite{MA-CN:18,MZ-CGC:08,CN-JC-GJP:15-acc}.
More specifically, we are interested in algorithms that use Voronoi partitions to solve the deployment problem. The authors in~\cite{CN-JC:11-auto} proposed a self-triggered Voronoi partitioning algorithm, the authors in~\cite{ZZ-XZ-LJ:18,NH-YM-ST:15} proposed an event-triggered Voronoi partitioning algorithm, and the authors in~\cite{CN-JC:15-crc} proposed a team-triggered Voronoi partitioning algorithm. The major drawback of these algorithms is the fact that when a trigger occurs, the agents must make a request for the information they require and other agents must respond to those requests i.e. instantaneous communication and synchronous actions by agents. { We would like to note that the instantaneous communication (immediate requests and immediate responce) can also be viewed in the context of sensing capabilities as well. This is due to the fact that both methods obtain information from neighbors instantaneously when needed. Thus, if the agents have sensing capabilities, the agents can sense the surrounding environment to localize other agents rather than communicating with them via message transmissions.} In this work, we propose a one-way communication model {that is fully asynchronous} where agents decide when to initiate communication by only broadcasting their information to others. This is similar to the communication model proposed in~\cite{GSS-DVD-KHJ:13} for multi-agent average consensus. The broadcasting communication model can be considered a more practical strategy given that it is an asynchronous model that reduces the amount of communication traffic seen over a SWS by eliminating the need for response messages to occur.
Additionally, since agents are restricted from requesting information or sending acknowledgments, they are unable to determine if their information has been sent to all intended recipients. This is in contrast to \cite{CN:17,CN-JC:11-auto} where agents communicate back-and-forth instantaneously to ensure sufficient information exchange. {Also, the works in \cite{ZZ-XZ-LJ:18,NH-YM-ST:15} assume that the agents are aware of their neighbors and can precisely determine the communication range required, which is not a practical assumption.}
Instead, we will design our solution to allow agents to independently determine the sufficient broadcasting range required to ensure that their information is shared with the intended recipients.

\paragraph*{Statement of contribution}
{
The main contribution of our work is the design of a fully asynchronous communication \algo~that allows agents to reach the set of locally optimal solutions to the deployment problem.
Unlike all existing triggered solutions, our event-triggered broadcasting algorithm relies on agents to decide when to broadcast (push) information to others in the network without the need for a response from other agents. As a result, our communication model relaxed the assumptions of instantaneous communication and synchronous actions by agents (request and response actions).
In addition, given the challenges of the communication constrains, we develop a distributed control algorithm such that the agents can determine the sufficient communication range required to reach all intended recipients.
Finally, our distributed algorithm is shown to achieve a level of performance that is similar to that of a deployment strategy that uses continuous or periodic communication methods.
}

\section{Preliminaries}
Let $\real$,~$\npreal$,~and~$\npint$ denote the sets of real, non-negative real, and non-negative integer numbers, respectively. Let $|\cdot|$ be the cardinality of a set. Also, we denote the Euclidean distance between two points $p,q \in \real^2$ by $\TwoNorm{p-q}$.

Let $Q$ be a convex polygon in $\real^2$ with a probability density function $\phi: Q \rightarrow \npreal$ that maps the probability of a spatial action or event occurring at point $q\in Q$. 
The \textit{mass} and \textit{center of mass} of $Q$ with respect to the density function $\phi$ are
\begin{align*}
M_Q= \int_Q \phi(q) dq \quad \text{and} \quad C_Q=\frac{1}{M_Q}\int_Q q\phi(q) dq,
\end{align*}
respectively. For a bounded set $Q \subset \real^2$, 
the \textit{circumcenter}, $cc(Q) \in \real^2$, is the center of the closed ball of a minimum radius contained in $Q$, and the \textit{circumradius}, $cr(Q) \in \npreal$, is the radius of the closed ball.
Then, let the closed ball centered at $q$ with a radius $r$ be $\OB(q,r)$. 
 
\subsection{Voronoi partition}
In this subsection, we briefly present some concepts necessary for the development of the \algo; further details on \textit{Voronoi partitions} can be found in~\cite{AK-BB-KS-SNC:09}.
For a convex polygon, $Q \subset \real^2$, 
let $P=\{\pone,\dots,\pN\}$ denote the locations of $N$ agents in $Q$ and let $\ID=\{1,\dots,N\}$ 
denote the set of identification numbers corresponding to the $N$ agents with locations $P$. The set $Q$ can be partitioned into $N$ polygons $\VV(P)=\{\voone,\dots,\voN\}$ such that the union of their disjoint interiors is $Q$. 
The Voronoi cell of agent~$i$ is formally defined as
\begin{align}\label{eq:voronoi}
 \voi=\{q \in Q \;| \;\TM{q-\ppi} \leqslant \TM{q-\pj}\; \forall\; i \neq j\} .
 \end{align}
When all agents are positioned at the centroids of their Voronoi cells, i.e., $\ppi=\cvi, \; \forall i \in \ID$, the agents' locations $P=(\pone,\dots,\pN)$ are said to be in a \textit{centroidal Voronoi configuration}.
Furthermore, when the intersection of two Voronoi cells $\voi$ and $\voj$ generated by the points $\ppi$ and $\pj$ is non-empty, $\voi \cap \voj \neq \emptyset$, the agents~$i$ and $j$ are Voronoi neighbors.
The set of Voronoi neighbors for the $i\text{th}$~agent is denoted by $\nei$.  

\subsection{Space partition with uncertain information}\label{se:uncertain}
In order to compute the Voronoi cells defined by~\eqref{eq:voronoi}, each agent requires the exact positions of their neighboring agents.
{If agents do not continuously communicate and/or sense their surroundings}, the exact locations of neighbors may not be available. Then, the agents must rely on inexact information to approximate their Voronoi cells. Here we discuss the concept of space partitioning when an agent's knowledge of its neighbors' positions is uncertain. Two methods for space partitioning under uncertain conditions are utilize, the \textit{guaranteed} and the \textit{dual-guaranteed Voronoi diagrams}~\cite{CN-JC:11-auto,CN-JC-GJP:15-acc,WE-JS:08,JM-MA-MM:09}. 

Let $\XX=\{\xone,\dots,\xN\} \subset Q$ be a
collection of compact sets containing the true positions of agents with~$\ppi \in \xi, \;\forall i\in \ID$. The set~$\xi$ is considered to be a region of uncertainty and represents all the possible points in $Q$ where agent~$i$ could potentially be located. If the location of agent~$i$ does not contain uncertainty, then the set~$\xi$ is simply a singleton with $\xi=\{\ppi\}$. 

\begin{figure}
\begin{center}
\subfigure[]{\includegraphics[width=4cm]{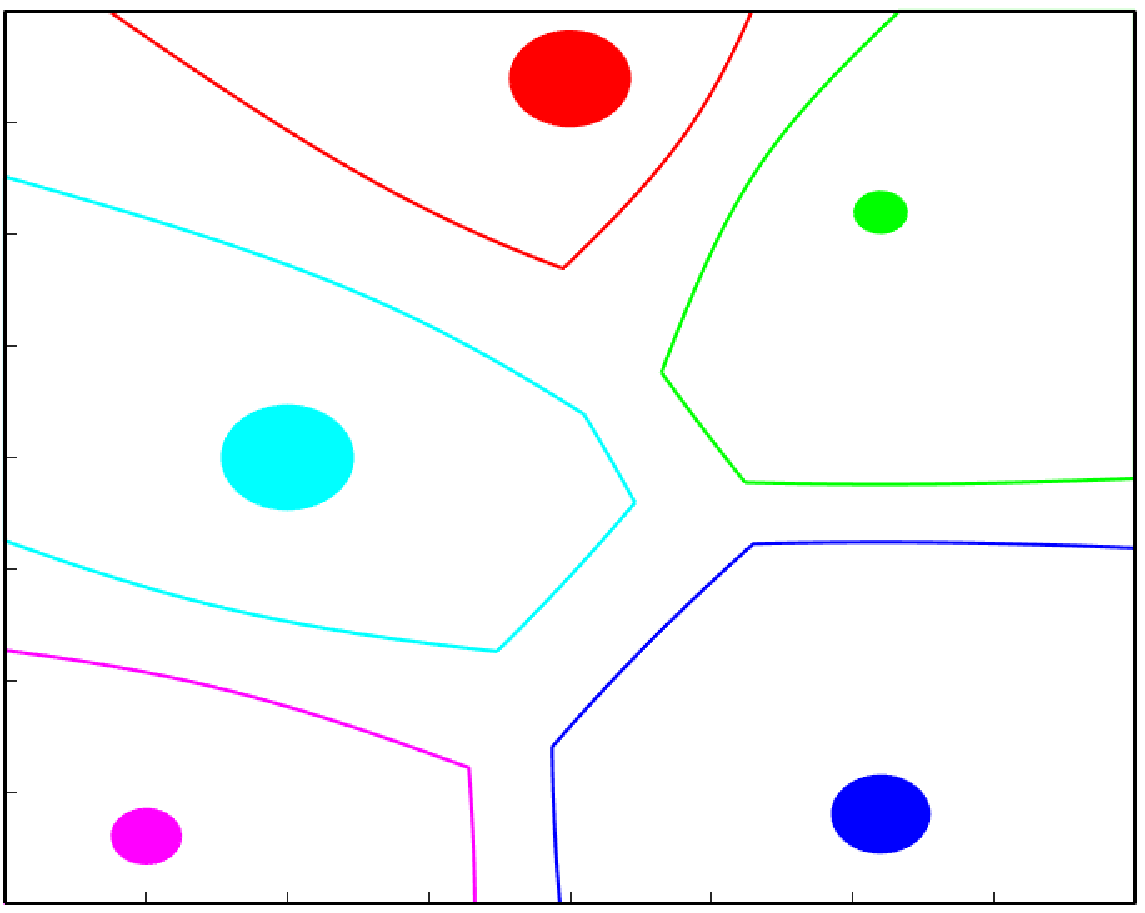}}
\subfigure[]{\includegraphics[width=4cm]{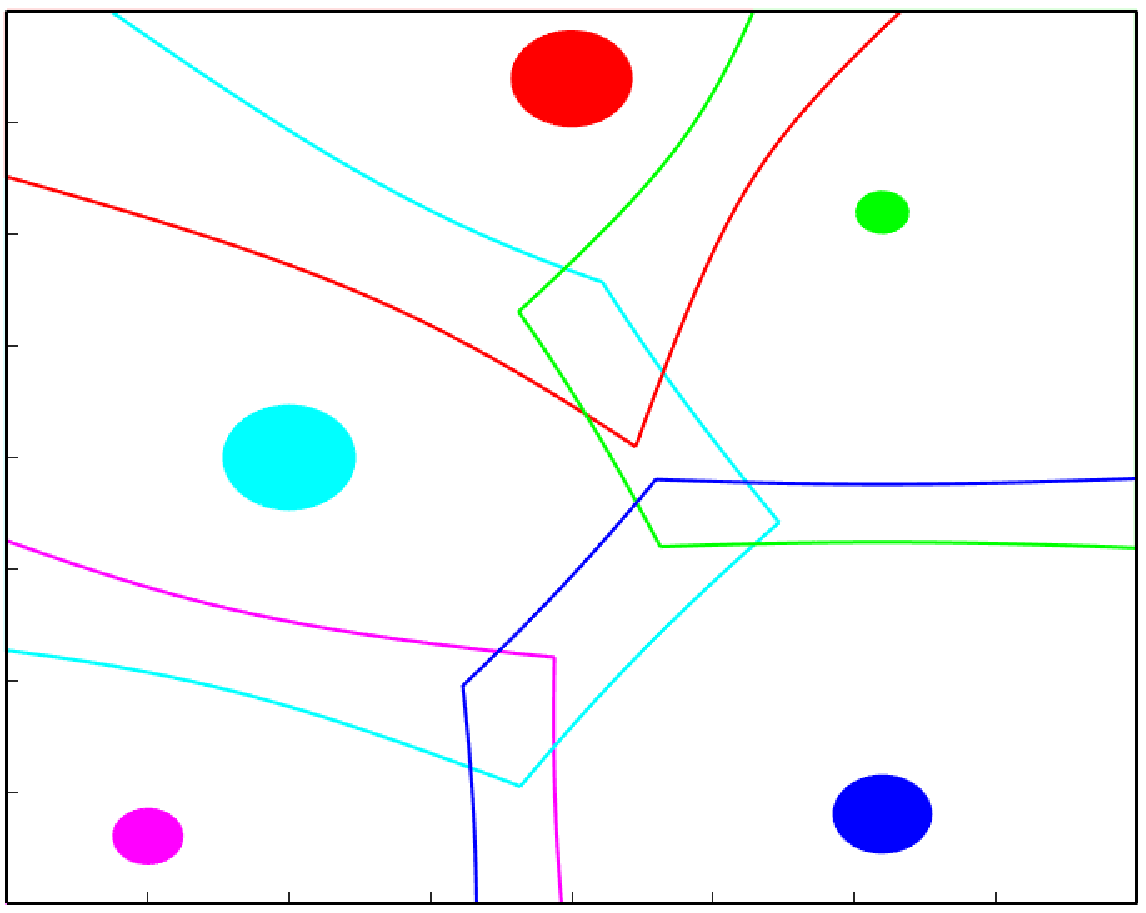}}
\caption{(a). Guaranteed and (b). dual-guaranteed Voronoi diagrams}\label{fig:1}
\end{center}
\end{figure}

The guaranteed Voronoi diagram, also known as the fuzzy Voronoi diagram~\cite{JM-MA-MM:09}, is the collection $g\VV(\XX) = \{\gvone,\dots,\gvn\}$ of guaranteed Voronoi cells generated by the uncertainty regions of~$\XX$. The guaranteed Voronoi cell~$\gvi$ for agent~$i$ is the set of points that are guaranteed to be closer to agent~$i$ than any other agents.
Formally,
\begin{align*}
\gvi=\{q \in Q \;|\; \underset{\xli\in \xi}{\max}\TM{q-\xli} \leq \underset{\xlj\in \xj}{\min}\TM{q-\xlj}\;\forall i\neq j\}.
\end{align*}
Note that in general, the guaranteed Voronoi cell is a subset of the Voronoi cell, i.e. $\gvi\subset \voi\;\forall i\in \ID$, 
which implies that the guaranteed Voronoi diagram is not a partition of $Q$.
Fig.~\ref{fig:1}.(a). provides an example of a guaranteed Voronoi diagram consisting of five agents with uncertainty regions represented by circles. 

Complementary to the guaranteed Voronoi diagram, the dual-guaranteed Voronoi diagram, first introduced in~\cite {CN-JC:11-auto}, is the collection $dg\VV(\XX)=\{ \dgvone,\dots ,\dgvn\}$ of dual-guaranteed Voronoi cells generated by the uncertainty regions of $\XX$. The dual-guaranteed cell $\dgvi$ for agent~$i$ is the collection of points such that any point outside the cell is guaranteed to be closer to all other agents than to agent~$i$. Formally,
\begin{align*}
\dgvi=\{q \in Q \;|\; \underset{\xli\in \xi}{\min}\TM{q-\xli} \leq \underset{\xlj\in \xj}{\max}\TM{q-\xlj}\;\forall i\neq j\}.
\end{align*}
Note that in general the Voronoi cell is a subset of the dual-guaranteed Voronoi cell, i.e.  $\voi \subset \dgvi \;\forall i\in \ID$. Fig.~\ref{fig:1}.(b). shows the dual-guaranteed Voronoi diagram of five agents and their uncertainty regions.

\subsection{Facility location}
In this subsection, the locational optimization function presented in \cite{JC-SM-TK-FB:02-tra} is discussed. 
The function that quantifies the sensing performance of an agent located at point~$\ppi$ to a point of interest $q \in Q$ is given by,
\begin{align}
f(\TM{q-\ppi}) = \TM{q-\ppi}^2\label{eq:f}.
\end{align}
This function measures the sensing quality of an agent based on its distance to a given point $q$. As the distance from an agent position $\ppi$ to the point of interest $q$ decreases, the sensing performance for agent~$i$ at point $q$ increases. For a density function $\phi(q): Q \rightarrow \npreal$ that captures the likelihood of an action occurring at $q$, the total network performance of $N$ agents at fixed positions $P$ is given by,
\begin{align}
\HH(P)=E_{\phi}\Big[ \underset{i\in \ID }{\min}\TM{q-\ppi}^2 \Big].
\end{align}
We assume $\phi(q)$ is provided to agents prior to deployment. The $\HH$ function is beneficial when the agents are the closest to the actions that they are responsible for, and it has been used in a number of applications previously, such as in event detection and resource allocation~\cite{QD-VF-MG:99,FB-JC-SM:09}. 
Given a Voronoi partition, each agent will be responsible for the points that are closer to itself than to any other agent, and now, the objective function with respect to a Voronoi partition as in \cite{JC-SM-TK-FB:02-tra} is written as
\begin{align} \label{eq:obj}
\HH(P)=\sum_{i=1}^{N}\int_{\voi} \TM{q-\ppi}^2 \phi(q) dq.
\end{align}
For a distributed control algorithm where the agents have local or limited information, the agents optimize the objective function by moving toward the centroid of their Voronoi cells. For $P'\in Q$ with $\TM{\ppi'-\cvi}\leq \TM{\ppi-\cvi}$ for all $i\in \ID$, then
\begin{align}\label{eq:pro1-3}
\HH(P',\VV(P))\leq \HH(P,\VV(P)).
\end{align}
{In addition, when agents are located at the centroids of their Voronoi cells, $(\pone,\dots,\pN)=(\cvone,\dots,\cvN)$, the objective function is considered to be in a locally optimized state \cite{JC-SM-TK-FB:02-tra}. In other words, when agents are located at the centroids of their Voronoi cells, the agents cannot improve the objective function by moving in any direction. This is commonly referred to as Nash equilibrium or a local minima result.}

\section{Problem Statement}
Let $P=\{\pone,\dots,\pN\} \in Q^N$ be the location of~$N$ agents moving in a convex polygon $Q \subset \real^2$. In this work we consider first-order dynamics for each agent,
\begin{align}\label{eq:dy}
\dot{p}_i=u_i,
\end{align}
where $u_i$ is the control input for agent~$i$, which is constrained by $\TwoNorm{u_i} \leq \vmax$ with $\vmax$ being the maximum speed for all agents.
{We assume that all agents have knowledge of the density function $\phi(q)$ before starting the deployment task. In such a case, the density function can be seen as a populated area or as a higher crime area in a city that require more surveillance and monitoring. Similar to \cite{CN-JC:11-auto,JC-SM-TK-FB:02-tra}, the $\phi(q)$ can be fixed during the deployment task. In addition, in case the density function is unknown, there have been many works to proper estimate a density function using, but not limited to, distributed Kriged Kalman filter and neural network as in \cite{JC:09,AD-SE-SMAS-MBM:12}. The authors in these works consider the scenario where the agents response to a dynamic change in an environment that is modeled by unknown dynamic density function.}
{ Also, we assume that agents are capable of maintaining speeds between $0$ and $s_{\max}$ for some duration of time. It is worth pointing out that many ground, under water, and unmanned aerial vehicles have this capability.}

Our objective is to achieve a locally optimal (local minimum) value of the objective function~$\HH$ in \eqref{eq:obj} with {a fully asynchronous communication model while also taking a more communication-aware approach.} 
More specifically, {rather than requiring instantaneous communication and/or synchronous actions by agents, as in many similar triggered deployment algorithms~\cite{CN-JC:11-auto,CN-JC:15-crc,ZZ-XZ-LJ:18,NH-YM-ST:15}}, 
we consider a {fully asynchronous} broadcast model without acknowledgment of or requests for information where the agents must actively choose a broadcast radius to transmit their messages similar to \cite{GSS-DVD-KHJ:13}. 
Additionally, we want to reduce the amount of communications across the SWS by having the agents only broadcast messages when necessary, rather than continuously or periodically. This means the agents will need to determine exactly when to broadcast messages to one another, with what distance, and how to move in the environment based on locally available information. 

Formally, let~$\{\tiell\}_{\ell \in \integernonnegative} \subset \realnonnegative$ be the sequence of times (to be determined on-line) at which agent~$i$ broadcasts its position~$\ppi$ to other agents in the network. The radius at which the message is broadcast is~$\riell > 0$, meaning any agent~$j$ that is within~$\riell$ of~$\ppi(\tiell)$ will receive the message.

Let~$\pij(t) \in Q$ be the last known position of agent~$j$ by agent~$i$ at any given time~$t \in \realnonnegative$. Note that if agent~$i$ has never received information from some agent~$j$, then~$\pij(t) = \emptyset$. Then, given a sequence of broadcast times~$\{\tjell\}_{\ell \in \integernonnegative}$ and broadcast radii~$\{\rjell\}_{\ell \in \integernonnegative}$ for each agent~$j \in \ID\setminus \{i\}$,
the memory of each agent~$i$ is updated according to
\begin{align}\label{eq:comm}
\pij(t) = \begin{cases} \pj(\tjell) & \text{if } \ppi(\tjell) \in \cball{\pj(\tjell)}{\rjell}, \\ \pij(\tjellm) & \text{otherwise}, \end{cases} 
\end{align}
for $t \in [\tjell, \tjellp)$. 

The goal is now to devise a distributed coordination strategy 
{such that the agents converge to a centroidal Voronoi configuration, which locally optimizes~$\HH$ in \eqref{eq:obj}, as in \cite{JC-SM-TK-FB:02-tra}. However, we do so with an asynchronous communication aware model.}
To this end we are interested in broadcasting messages as minimally as possible, both in frequency and space, by only broadcasting when necessary. 

\begin{problem}\label{prob:1}
Given the dynamics~\eqref{eq:dy} and the communication model~\eqref{eq:comm}, find a distributed communication and control strategy to find a sequence of broadcasting times~$\{\tiell\}_{\ell \in \integernonnegative}$, broadcasting ranges~$\{\riell\}_{\ell \in \integernonnegative}$, and a control strategy~$u_i(t)$ for all~$i \in \ID$ that {drives the agents to their Voronoi centroids in order to locally} optimize the objective function~$\HH$ in \eqref{eq:obj}.
\end{problem}

\section{Event-Trigger Algorithm Design}
When an agent has knowledge of the exact locations of its Voronoi neighbors, the agent can compute its exact Voronoi cell and the exact location of the cell's centroid. This allows the agent to move directly toward the Voronoi cell centroid. This will monotonically optimize the objective function~$\HH$ in \eqref{eq:obj} as in~\cite{JC-SM-TK-FB:02-tra}. Unfortunately with a {fully asynchronous event-triggered} broadcasting communication model, the agents do not know when they will receive new information from other agents, nor can they request the information when they require it as is done in~\cite{CN-JC:11-auto,CN-JC:15-crc,ZZ-XZ-LJ:18,NH-YM-ST:15}. Instead, they must rely on the information that they possess at each moment in time in order to determine exactly how to move and when to initiate data transfers to others. Although, it is known from~\cite{CN-JC:11-auto} that agents only need information from their Voronoi neighbors {to move to their Voronoi centroids}, the major challenge now is that the agents do not necessarily know their actual Voronoi neighbors. This is due to the continuously increasing uncertainty that exists with respect to the positions of other agents when not communicating. This is combined with the fact that agents cannot gain knowledge of other agent positions by requesting the information.

\subsection{Information Required and Memory Structures} \label{se:MBR}
Since agents do not exchange information with each other on a continuous basis, each agent must maintain the most current state information they have received from other agents, as well as a method to model the uncertainty that evolves over time. The data structure that allows the agents to compute the uncertainty regions of any other agent using the most recently received information with respect to agent $j \in \ID \setminus\{i\}$ is the following. One, the time instance $\tij$ when agent~$i$ last successfully received information from agent~$j$. Two, the position $\pij=\pj(\tij)$ of agent~$j$ received at time $\tij$, and three, the speed promise~$\vij=\vj(\tij)$ by agent~$j$ received at time $\tij$. The notion of the promise is borrowed from~\cite{CN-JC:16-tac}.
The speed promise~$\vij$ made to agent~$i$ by agent~$j$ states that agent~$j$ promises not to exceed the speed given by the promise~$\vij$. In other words, agent~$j$'s dynamics will follow~$\TwoNorm{u_j(t)} \leq \vij$, and this will hold for agent~$j$ until it broadcasts again. For now, we consider agent~$j$ only setting its speed to one of two values~$\vj \in \{0, \vmax \}$. This can be seen as two modes of operation signaling either `active' i.e. agent~$j$ is moving or `inactive' i.e. agent~$j$ is not moving by holding its current position. {This assumption helps simplify the convergence result in~Section~\ref{Se:ConAnalysis} and will be relaxed in Section~\ref{se:ex:cv} where agents can modify their speed and promises to any~$s_j \in [0, \vmax]$.}

We define $\dij = (\tij, \pij, \vij ) \in (\npreal\times Q \times \{0,\vmax\})$ as the information that agent~$i$ last received about any given agent~$j$, and $\dij = (\emptyset)$ if agent~$i$ has not received information from agent~$j$. This information allows agent~$i$ to construct a closed ball that guarantees to contain the $j$th agents' real location. For any time $t\geq \tij$, agent~$i$ knows that the $j$th agent did not move farther than $\vij(t-\tij)$ away from~$\pij$. Given the promise and the dynamics~\eqref{eq:dy}, the uncertainty region of the $j$th agent with respect to the information agent~$i$ has is formally defined as
\begin{align}\label{eq:unc}
\xij(t)=\OB (\pij,\vij(t-\tij)) .
\end{align}
We denote by $\dii = (\tii,\pii,\vii)$ the $i$th agent's information at its latest broadcast time $\tii$ and by $\di = (t,\ppi,\vi)$ the $i$th agent's current information. 
The $i$th agent's full memory at any given time is collected in
\begin{align}
&\ddi=(\dione,\dots,\diN)\in (\npreal\times Q \times \{0,\vmax\})^N \nonumber .
\end{align}
Additionally, the memory for all agents in the entire network is defined by
\begin{align}
\DD=\{\ddone,\dots,\ddN\}\in (\npreal\times Q \times \{0,\vmax\})^{N^2} .
\end{align}
With the memory structure defined, we are now ready to begin solving Problem~\ref{prob:1}. We begin by determining exactly which agents in the network should a particular agent {broadcast to} in order to complete the deployment task. Since it is known from~\cite{CN-JC:11-auto} that an agent's knowledge of its Voronoi neighbors is useful in solving a similar self-triggered deployment problem, we determine a method for agents to keep track of their Voronoi neighbors for the event-triggered broadcasting problem given the uncertainty regions. For a given agent~$i$ and its memory~$\ddi$, we propose the notion of \emph{dual-guaranteed neighbors} next.
Then, as long as agent~$i$ maintains some type of communication with its dual-guaranteed neighbors~$j \in \dgi$, we conclude that the agents will have a sufficient amount of information to allow them to {move toward their Voronoi centroids effectively locally} optimizing the objective function~$\HH$ in \eqref{eq:obj}. 

\begin{definition}[Dual-Guaranteed Neighbors]\label{def:DguaNei} 
{\rm
Given a set of uncertainty regions~$\XX = (\xone, \dots, \xN)$ such that~$\pj \in \xj$ for all~$j \in \ID$, a \emph{dual-guaranteed neighbor} of agent~$i$ is any agent~$j$ that can be made a Voronoi neighbor for at least one configuration of positions~$P \subset \XX$. Formally, the set of dual-guaranteed neighbors of agent~$i$ is
\begin{align*}
\dgi=\{j\in\ID \;|\; \exists P \subset \XX \text{ s.t. } j\in\nei   \} .
\end{align*} 
}
\end{definition}	

The dual-guaranteed neighbors are defined such that~$\nei \subset \dgi$ which is formalized in Lemmas~\ref{lm:ndg} below. 

\begin{lemma}[Dual-Guaranteed Neighbors] \label{lm:ndg}
Given a set of uncertainty regions~$\XX = (\xone, \dots, \xN)$ such that~$\pj \in \xj$ for all~$j \in \ID\setminus\{i\}$, if
\begin{align*}
 \dgvi(\XX) \cap \dgvj(\XX)  \neq \emptyset ,
\end{align*}
then $\exists P \subset \XX$ such that agent~$j\in\dgi$ can be a Voronoi neighbor of agent~$i$.
\end{lemma}

\begin{IEEEproof}  
In appendix~\ref{app:A}.
\end{IEEEproof}

Moreover, an agent~$i$ might not have any information about some agents in the network if they have not communicated yet. This is troublesome since the computation of $\gvi$ and $\dgvi$ are based on the availability of the uncertainty regions that guarantee~$\pj \in \xij$ for all times. An even bigger challenge given the setup of our problem is the possibility that some agents are in communication for some period of time before stopping altogether, since it is no longer necessary. 
In case that agents do not communicate with all other agents, we need to determine additional conditions to ensure that, with partial information from $\ddi$, 
$\nei \subset \dgi$ is guaranteed.
To address this, we define a map~$\pi_{\JJ} : (\npreal\times Q \times \{0,\vmax\})^{N} \rightarrow (\npreal\times Q \times \{0,\vmax\})^{|\JJ|}$ for any~$\JJ \subset \ID$ that extracts information corresponding only to agents~$j \in \JJ$ from~$\ddi$. Formally,
\begin{align*}
\pi_{\JJ}(\ddi) = \cup_{j \in \JJ} \{ \dij \} .
\end{align*}

Furthermore, we require $\pi_\JJ(\ddi)$ to be sufficient to guarantee
$\nei \subset \dgi$.
Intuitively, this means that all $k\not\in\JJ$~agents should be sufficiently far away from agent~$i$ that agent $k$ is not a Voronoi neighbor of agent $i$ given it $\pk(t)$. This is formalized in Corollary~\ref{coro} below. 

\begin{corollary}[Condition on Voronoi Neighbors Sets] \label{coro} 
Given a set of uncertainty regions~$\XX = (\xi\cup \{\xij\}_{j\in\JJ})$ such that~$\xi=\overline{B}(\ppi,0)$ and $\pj \in \xij$ for all~$j \in \JJ$, if
\begin{align*}
\pk \not\in \overline{B}\left(\ppi,2\cdot\underset{q\in\dgvi(\XX)}{\max}(\TM{q-\ppi})\right)\;\forall k\in\ID\setminus\JJ \setminus\{i\},
\end{align*}
then $\nei\subset\dgi$ is guaranteed.
\end{corollary}

\begin{IEEEproof}  
In appendix~\ref{app:B}.
\end{IEEEproof}

It is important to emphasize that for any $j\in\dgi$, if agent~$j$ is a dual-guaranteed neighbor of agent~$i$, it does not imply that agent~$i$ is a dual-guaranteed neighbor of agent~$j$ since $\ddi \neq \ddj$.
{It is worth noting that agents do not need to have knowledge of the number of agents in the network. In order to describe the algorithm in a simpler manner, $\ddi$ contains a place holder for each agent in the network (even if the actual execution of the algorithm will never need to access all of this information), and as further explained, the agents only need information from the dual-guaranteed neighbors.}
For convenience, we let 
\begin{align*}
\gvij(\JJ)&=\gvj(\{\xij \}_{j \in \JJ}\cup \xi),\\
\dgvij(\JJ)&=\gvj(\{\xij \}_{j \in \JJ}\cup \xi).
\end{align*} 
The objects in this subsection are summarized in Table~\ref{tb:para} for agent~$i,j\in \ID$.

\begin{centering}
\begin{table}[h]
  \setlength{\tabcolsep}{3pt}
  \begin{tabular}{|p{60pt}|p{175pt}|}
\hline
$\ppi\in \real^2$ & agent~$i$'s location\\
$\pij\in \real^2$  &  agent~$j$ broadcasted location to agent $i$\\
$\vi\in \{0,\vmax\}$  &  agent~$i$'s speed\\
$\vij\in \{0,\vmax\}$  &  agent~$j$ promised speed to agent~$i$\\
$\tij\in \npreal$  &  agent~$j$ broadcasting time to agent~$i$\\
$\nei\subset \ID$ &  agent~$i$'s Voronoi neighbors\\
$\dgi\subset \ID$&  agent~$i$'s dual-guaranteed neighbors\\
$\ddi\in \DD$ &  agent~$i$'s full memory\\
$\dij\in \ddi$ &  agent~$j$'s information in the $i$th agent's memory\\
$\xij\subset Q$ &  agent~$j$'s uncertainty region given $\dij$ \\
$\gvij(\JJ) \subset Q$&  agent~$j$'s guaranteed Voronoi cell with respect to\\ &  $\pi_\JJ(\ddi)\cup \di$ information\\
$\dgvij(\JJ) \subset Q$&  agent~$j$'s dual-guaranteed Voronoi cell with respect to\\ &  $\pi_\JJ(\ddi)\cup \di$ information\\
\hline
 \end{tabular}\caption{Agent~$i$ model definitions}\label{tb:para}
\end{table}
\end{centering}

\subsection{Motion Control Law} \label{se:MC}
The motion control law defines a method to generate trajectories for the agents that allows them to contribute positively to the deployment task.
{ This is accomplished by having agents move towards the midpoint between the centroids of their guaranteed and dual-guaranteed Voronoi cells.} 
Due to uncertainties, agents must use the information pertaining to their dual-guaranteed neighbors, which is guaranteed to contain all their Voronoi neighbors. Since $\nei\subset\dgi$ by Lemma~\ref{lm:ndg}, the agents guarantee computing the guaranteed and dual-guaranteed Voronoi cells such that $\gvii(\dgi) \subset \voi \subset \dgvii(\dgi)$. 
Now let us informally describe the idea behind it here.

At each time-step, agent~$i$ uses partial information, its dual-guaranteed neighbors' information, to determine if it can move in a manner to optimize~$\HH$ in \eqref{eq:obj}. 
{ If so, it computes the centroids of the guaranteed and dual-guaranteed Voronoi cell and then moves toward the midpoint between them.} 
Otherwise, it does not move and waits until it receives sufficient information to initiate the continuation of motion. 

{In general, it is preferable that agents move toward their Voronoi centroids directly. However, this requires perfect information at all times. Instead, we establish a convex set~$\CC_i(\dgi)$ that is guaranteed to contain the true centroid~$C_{V_i}$ based on the information~$\dgi$ available to agent~$i$. 
This set can then be used not only to determine how to move but also exactly when updated information is needed. This is formalized next in Proposition~\ref{pro:newBound}.

\begin{proposition}\label{pro:newBound}
Given~$D_j^i = (t_j^i,p_j^i,s_j^i)$ for all~$j \in \dgi$, let
\begin{align}\label{eq:bound1}
\CC_i(\dgi) = \cball{\cdgi}{\bnd} \cap \cball{\cddgi}{\bnd},
\end{align}
where
\begin{align*}
\bnd=2cr_{\dgvii(\dgi)} \left(1- \frac{M_{\gvii(\dgi)}}{M_{\dgvii(\dgi)}} \right),
\end{align*}
then~$C_{V_i} \in \CC_i(\dgi)$. 
\end{proposition}

Proposition~\ref{pro:newBound} says that although based on uncertain information the exact location of~$C_{V_i}$ cannot be determined, its distance from the centroids of the guaranteed~$\cdgi$ and dual guaranteed~$\cddgi$ can be upper-bounded by the same quantity~$\bnd$. The set~$\CC_i$ is then just the intersection of the two balls centered at~$\cdgi$ and~$\cddgi$ with radii~$\bnd$. It is then easy to see that the set~$\CC_i$ is convex, and thus while agent~$i$ is outside of this set, it can guarantee to be getting closer to the true centroid~$C_{V_i}$ by simplify moving towards~$\CC_i$. 
With that being said, we define $m_i$ as the closest point on $\CC_i(\dgi)$ to agent $i$. Formally,
\begin{align}\label{eq:medpoint}
m_i=\underset{q\in \CC_i(\dgi)}{\arg\min} \|q-\ppi\|,
\end{align}
and as long as 
\begin{align}\label{Condition}
\ppi\neq m_i,
\end{align}
the agent has good information to move. 
Then, our motion control becomes
\begin{align}\label{eq:ui}
u_i= \vi \frac{\mpi - \ppi}{\TM{\mpi - \ppi}},\; \forall \; i\in\ID,
\end{align}
where
\begin{align*}
\vi=\begin{cases} 
\vmax \hspace*{10 mm} \text{if condition \eqref{Condition} is true} , \\
0 \hspace*{15.3 mm}\text{otherwise }.
\end{cases}
\end{align*}
}

\begin{figure}
\vspace*{-6ex}
\vskip3em
{\includegraphics[trim={0cm 1.4cm 0 1.2cm},clip,width=10cm]{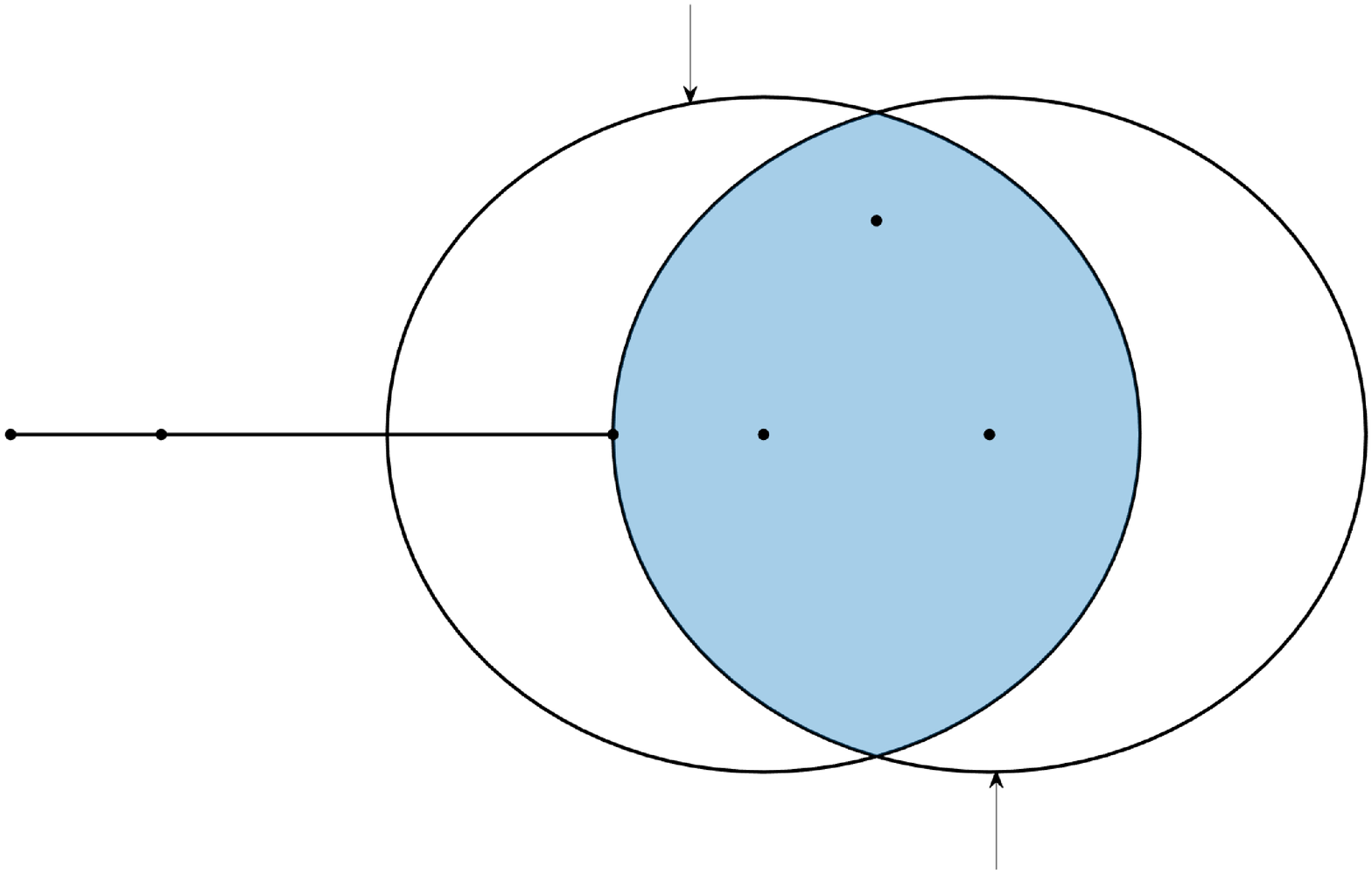}}
\put(-240,55){\small $\ppi$}
\put(-220,55){\small $\ppi$}
\put(-135,52){\small $C_1$}
\put(-105,52){\small $C_2$}
\put(-125,85){\small $C_{V_i}$}
\put(-175,125){\small $\cball{\cdgi}{\bnd}$}
\put(-130,-5){\small $\cball{\cddgi}{\bnd}$}

\centering
\caption{Example of a point $\ppi'$ that an agent can move to given $C_1=\cdgi$, $C_2=\cddgi$ and $\bnd$}\label{fig:2}
\end{figure}

\vspace*{-0ex}
\begin{algorithm}[htb]
 {\footnotesize At any time $t>0$, agent~$i\in  \ID$ performs:
  \null\hfill\null \vspace*{-0ex}
  \begin{algorithmic}[1]
   \STATE sets $D =\pi_{\dgi}(\ddi)\cup \di$
   \STATE computes $\XX(D)$ as in \eqref{eq:unc}
   \STATE computes $L=\gvi(\XX)$ and $C_L$
   \STATE computes $U=\dgvi(\XX)$ and $C_U$
   \STATE computes $\mpi$ as in \eqref{eq:medpoint}
   \STATE computes $u_i$ as in \eqref{eq:ui}
  \end{algorithmic}
 }
 \caption{\hspace*{-.5ex}: \small Motion Control Law}\label{tab:mcl}
\end{algorithm}\vspace*{-0ex}

\subsection{Decision Control Law} \label{se:DC}
Equipped with a motion control law, we now need to design a communication protocol that provides sufficient information to the motion control law to properly do its job.
This is nontrivial because the agents do not have control over when they will receive updated information from others, but instead can only choose when they send information. {
Thus, rather than designing the decision control law for agent~$i$ based on when \emph{it} needs information, our decision control law would ideally be designing in terms of when its neighbors need information. However, since agents in general may not know exactly what their neighbors need (due to distributed information), we instead propose a law that broadcasts messages as minimally as possible while still ensuring the entire network converges to the desired set of states. 
}

Starting with the case when agent~$i$ increases its speed from $0$ to $\vmax$ without broadcasting, the neighbors $j\in\dgi$ will not be able to capture the correct uncertainty region of agent~$i$.
This can lead to failure of achieving the deployment task if the case is not handled appropriately. To handle this situation, the dual-guaranteed neighbors must be informed about the new change in speed as soon as it increases. This makes it possible for other agents to appropriately manage the uncertainty that they possess for agent~$i$'s location by modifying the rate of change at which the uncertainty evolves. 
Therefore, agents shall broadcast their information when their speed increases i.e. change from $0$ to $\vmax$. 

For the case when agent~$i$ changes its speed to $0$ without broadcasting, the uncertainty about its location held by other agents will increase. As a result, its uncertainty will eventually become larger than necessary. Consequently, the dual-guaranteed neighbors will no longer have enough reliable information to continue to move. To prevent this scenario, the $i\text{th}$ agent's dual-guaranteed neighbors need to know when the agent's speed has been set to $0$. This allows the neighbors to halt the expansion of the uncertainty region for agent~$i$. Therefore, agents shall broadcast their state information when they change their speeds to $0$. 

For the case when agent~$i$ realizes that it has a new dual-guaranteed neighbor $j\in\dgi$, both agents~$i$ and~$j$ must know about each other because they can be Voronoi neighbors.  
Thus, agent~$i$ must broadcast as soon as it gets information from the new dual-guaranteed neighbor to ensure agent $j$ has its information. Furthermore, we would like to note that for a deployment problem, as the agents move away from each other, the set $\dgi$ gets smaller, and if this case happens, it only occurs finite times during the deployment task.

Algorithm~\ref{tab:dcl} formalizes the decision control law.
\vspace*{-0ex}
\begin{algorithm}[htb]
 {\footnotesize At any time $t>0$, agent~$i\in  \ID$ performs:
  \null\hfill\null \vspace*{-0ex}
  \begin{algorithmic}[1]
   \IF{$\vi > \vii$}
   \STATE broadcasts $\di$
   \ELSIF{$\vi=0$}
   \STATE broadcasts $\di$
   \ELSIF{agent~$i$ has a new dual-guaranteed neighbor}
   \STATE broadcasts $\di$
   \ENDIF
  \end{algorithmic}
 }
 \caption{\hspace*{-.5ex}: \small Decision Control Law}\label{tab:dcl}
\end{algorithm}\vspace*{-0ex}

\subsection{Broadcasting Range}
Equipped with a motion control law and a method for determining exactly when broadcasting new information is necessary, we are interested in determining the minimum broadcasting range required. Since $\nei \subset \dgi$ by Lemma~\ref{lm:ndg}, we would like to broadcast to all agents $j\in\dgi$. This guarantees to send the information to all agents $j\in\nei$.

Thus, we aim to find the minimum radius for an agent to broadcast in order to guarantee that all $j \in \dgi$ are reached. 
This can be achieved by finding the distance from an agent~$i$ at $\ppi$ to its farthest dual-guaranteed neighbor. By Corollary \ref{coro}, any agent~$j$ s.t. 
\begin{align*}
\pj \not\in \overline{B}\left(\ppi,2\cdot\underset{q\in \dgvii(\dgi)}{\max}(\TM{q-\ppi})\right)
\end{align*}
is guaranteed not to be a dual-guaranteed neighbor. Therefore, the agent can broadcast with a distance of $(2\cdot\underset{q\in \dgvii(\dgi)}{\max}(\TM{q-\ppi})),$ which is sufficient to transmit information to all $j \in \dgi$. By Lemma~\ref{lm:ndg}, this guarantees that all the Voronoi neighbors of agent~$i$ receive agent~$i$'s information while agent~$i$ only uses the information provided by the data in $\dgvii(\dgi)$ to achieve this result.
The computation of the broadcasting range by agent~$i$ is formally defined as
\begin{align}\label{eq:rad}
&\ri(\dgi)=2\cdot\underset{q\in \dgvii(\dgi)}{\max}(\TM{q-\ppi(\tiell)}).
\end{align}

\section{Event Triggered Broadcasting Algorithm}\label{se:algo}
In this section, we combine the motion control law, decision control law and broadcasting range assignment to synthesize the {fully asynchronous communication }\algo. The goal is for the agents to converge to to the set of centroidal Voronoi configurations in order to locally optimize $\HH(P)$ given the communication model presented thus far. 

Given the nature of our problem, agent $i$ and $j$ may communicate for some duration of time and then stop when it is no longer necessary. However, the sets $X_i^j$ and $X_j^i$ will continue to expand even if communication is discontinued that effect the motion of both agents. To address this issue, we introduce a new set of neighbors called the potential neighbor set $\pni$ of agent~$i$ that exclude the dual-guaranteed neighbors of the $i$th agents which are guaranteed not to be Voronoi neighbors due to discontinuation of communication. In addition, the potential neighbor set satisfies $\NN_i\subset\pni\subset\dgi$. 
Now, let us introduce the mechanism to update the potential neighbor set when an agent receives new information as following: 
\begin{align}\label{eq:updateJ}
\pni^+=&\{ \{k\}_{k\in \pni'}\;|\; \dgvi(\XX)\cap \dgvk(\XX)\neq \emptyset\},
\end{align} 
where $\XX=\xi\cup\{\xij\}_{j\in\pni}$ and $\pni'=(\pni\cup j)$ if $j\not\in\pni$ and $\XX=\xii\cup\{\xij\}_{j\in\pni}$ and $\pni'=\pni$ otherwise.
Note, if an agent does not receive new information, $\pni$ does not change.

\begin{lemma}[Potential Neighbors Properties] \label{lm:J}
Under the decision control law, broadcasting range assignment and updating potential neighbors mechanism in \eqref{eq:updateJ}, $\nei\subset\pni$ is guaranteed at every event-time $\tiell$.
\end{lemma}

\begin{IEEEproof}  
In appendix~\ref{app:C}.
\end{IEEEproof}

\begin{proposition}\label{pro:algo}
Under the \algo, $\gvii(\pni) \subset \voi \subset \dgvii(\pni)$ is guaranteed at all times.
\end{proposition}

\begin{IEEEproof}  
In appendix~\ref{app:D}.
\end{IEEEproof}

Thus far, by Proposition \ref{pro:algo}, the algorithm ensures that the requirement for the motion control law is satisfied, { and the analysis of the asymptotic convergence properties is provided in the following section.}

Now, we proceed with an informal description of the proposed algorithm. Note that it is assumed that each agent knows their potential neighbors at time $t=0$ and consequently $\nei\subset\pni$. {Also, we assume that the agents are provided with the density function $\phi(q)$.}
Let us start with when the agent's speed is $0$. This implies agent~$i$ cannot contribute positively to the task. Therefore, it waits until it has a sufficient amount of information such that condition \eqref{Condition} holds where the uncertainties of the potential neighbors are computed using \eqref{eq:unc}. When agent~$i$ receives enough information to move, it sets its speed to $\vmax$ and broadcasts its state information at a distance $\ri(\pni)$ from its current position $\ppi$.

Furthermore, when agent~$i$ is able to contribute positively i.e. $\vi=\vmax$, it expands the uncertainties of its potential neighbors by the maximum rate of change, as in \eqref{eq:unc1}, and hold on the new received information until it broadcasts again. Then, the agent follows the motion control law \eqref{eq:ui} until condition \eqref{Condition} is not satisfied. We would like to note that an agent holds on the new information when it is in motion to allow it-self to follow the motion control law without the need to information from any agent $j\not\in\pni$ as explained in the proof of Proposition\ref{pro:algo}.  
Then, when condition \eqref{Condition} becomes invalid, the agent computes the uncertainties of its potential neighbors using \eqref{eq:unc} and checks for the condition \eqref{Condition}. If it is valid, the agent broadcasts its state information to a distance of $\ri(\pni)$ to ensure $\nei\subset\pni$ and repeats this until condition \eqref{Condition} is no longer valid. 
Next, the agent sets its speed to zero and broadcasts its state information to a distance of $\ri(\pni)$. Then, the agent waits until it gets sufficient information such that condition \eqref{Condition} is valid.
\begin{align}\label{eq:unc1}
\olxij(t)=\OB \left(\pij,\left(\vmax (t-\tii)+
\vij (\tii -\tij)\right)\right).
\end{align}
Algorithm~\ref{tab:ecl} formally describes the event-triggered broadcasting control law. 

\vspace*{-0ex}
\begin{algorithm}[htb]
 {\footnotesize Initialization at time $t=0$, agent~$i\in \ID$ performs:
  \null\hfill\null \vspace*{-2ex}
  \begin{algorithmic}[1]
  \STATE sets $D =\pid\cup D_i$
  \end{algorithmic}
}
  {\footnotesize At any time $t>0$, agent~$i\in \ID$ performs:
  \null\hfill\null \vspace*{-0ex}
  \begin{algorithmic}[1]
  \STATE updates $\pni$ using \eqref{eq:updateJ} 

  \IF{$\vi= 0$} 
  \STATE update $D =\pid\cup D_i$
  \STATE computes $\XX(D)$ as in \eqref{eq:unc}
  \ELSE
  \STATE computes $\overline{\XX}(D)$ as in \eqref{eq:unc1}
  \ENDIF
  \IF{$\vi = 0 \textbf{ and }$condition \eqref{Condition} is valid} 
  \STATE sets $\vi=\vmax$ 
  \STATE broadcasts $\di$ using \eqref{eq:rad} $\ri(\JJ_i)$ distance away 
  \ELSIF{$\vi\neq 0 \textbf{ and }$condition \eqref{Condition} is invalid} 
  
  \STATE update $D =\pid\cup D_i$  
  \STATE computes $\XX(D)$ as in \eqref{eq:unc}
  
  \IF{condition \eqref{Condition} is valid} 
  \STATE broadcasts $\di$ using \eqref{eq:rad} $\ri(\JJ_i)$ distance away 
  \ELSE
  \STATE sets $\vi=0$ 
  \STATE broadcasts $\di$ using \eqref{eq:rad} $\ri(\JJ_i)$ distance away 
  \ENDIF
  \ENDIF

  \IF{$t\neq t_i^i$ \textbf{and} agent~$i$ has new potential neighbor} 
  \STATE broadcasts $\di$ using \eqref{eq:rad} $\ri(\JJ_i)$ distance away 
  \ENDIF

  \STATE compute $u_i$ as in \eqref{eq:ui}
  \end{algorithmic}
 }
 \caption{\hspace*{-.5ex}: \small \algo}\label{tab:ecl}
\end{algorithm}\vspace*{-0ex}

\section{Convergence Analysis of the \algo}\label{Se:ConAnalysis}
{
In this section, we analyze the asymptotic convergence properties of the \algo.~Recall that our objective is to drive the agents to their Voronoi centroids because if the agents converges to the set of centroidal Voronoi configurations, the agents are locally optimized $\HH(P)$ \cite{JC-SM-TK-FB:02-tra}, which also means reaching a local optimal with respect to $\HH(P)$, where 
\begin{align}
\HH(P)=\sum_{i=1}^{N}\int_{\voi} \TM{q-\ppi}^2 \phi(q) dq.
\end{align}

\begin{proposition}
The agents' location evolving under the \algo~from any initial location configuration in $Q^N$ converges to the set of centroidal Voronoi configurations
\end{proposition}

\begin{IEEEproof}  
We know from \cite{JC-SM-TK-FB:02-tra} that 
\begin{align}
\dot{\HH}(P)=\sum_{i=1}^N \left(2M_{V_i}(\ppi-\cvi)\dot{\ppi}\right),
\end{align}
and now we want to show that $\dot{\HH}\leq 0$ under \algo
\begin{align*}
\dot{\HH}(P)&=\sum_{i=1}^N \left(2M_{V_i} (\ppi-\cvi)\frac{-\vi}{\TwoNorm{\ppi-\mpi}}(\ppi-\mpi)\right),\\
&= \sum_{i=1}^N \left(\frac{-2M_{V_i} \vi}{\TwoNorm{\ppi-\mpi}}(\ppi-\cvi)\cdot( \ppi-\mpi)\right).
\end{align*}
Under the motion control law if $\ppi\neq m_i$, the $(\ppi-\cvi)\cdot( \ppi-\mpi) \geq 0$. This implies $\dot{\HH}(P)<0$. In fact, $\dot{\HH}(P)$ is strictly negative if there is at least one agent moving toward its Voronoi centroid \cite{JC-SM-TK-FB:02-tra}. In case all agents are stationary such that $\vi =0\;\forall\; i\in\ID$, $\dot{\HH}(P)=0$. In addition, if the agents are stationary and stay stationary, it means that the agents know the exact location of their neighbors since $\vj =0\;\forall\; j\in\pni$, and it means that $\bnd=0$, $\cball{\cdgi}{\bnd}=\cdgi$, and $\cball{\cddgi}{\bnd}=\cddgi$ $\forall\;i\in\ID$. If condition \ref{Condition} is invalid for all agents, due to $\TwoNorm{\ppi-\mpi}=0$, this implies $\forall i\in\ID,\; \ppi=\cvi=\mpi=\cdgi=\cddgi$. Since $\ppi=\cvi\;\forall\;i\in\ID$, the agents converged to the set of centroidal Voronoi configurations. This conclude the proof.
\end{IEEEproof}
}

\section{Varying Speed Extensions}\label{se:ex:cv}
In this section, we consider an extension to the \algo~where agents can adjust their speeds to values other than zero or maximum. In other words we consider the case where agent~$i$ may take on a speed $\vi \in [0,\vmax]$. 
The main advantage of adjusting the speed of agents is to reduce the amount of communication between agents when an agent approaches its cell centroid while not effecting the asymptotic convergence. As the agents get closer to their Voronoi centroids, the agents start to communicate more frequently when they are moving at $\vmax$. Rather than allowing the agents to move with $\vmax$ all the time, an alternative approach is for agents to determine the appropriate speed to travel as they become close to their Voronoi centroids. 
For a variable speed adjustment approach, an agent's speed at any given instance in time can be described by the following,
\begin{align}\label{eq:vel}
\vi= \beta_i \vmax,
\end{align}
where $0<\beta_i \leq 1$ is determined online by the agents. Each agent updates $\beta_i$ based on a design parameter $\Delta_{TB}$, where $\Delta_{TB}$ represents a target time duration for agents to attempt to maintain movement prior to switching states and broadcasting information. The parameter $\Delta_{TB}$ is chosen prior to the start of the deployment task and the value chosen can depend on the anticipated communication traffic that will occur over the wireless network.
If agent~$i$ finds itself in a scenario where it stays in a moving state for less time than $\Delta_{TB}$, then this would imply that agent~$i$ is moving with a greater speed than desired. Thus, agent~$i$ will decrease it's speed by decreasing $\beta_i$. On the other hand, if the agent stayed in moving state greater than $\Delta_{TB}$, this would imply that agent~$i$ is moving with a speed that is slower than desired. Thus, agent~$i$ will increase its speed by increasing $\beta_i$. 
It is important to note that when agent~$i$ needs to increase the value of $\beta_i$, it must broadcast the new speed to its neighbors so that they can correctly capture the uncertainty region associated with agent~$i$ as explained in the first case in Section \ref{se:DC}.

In addition, by using the maximum speed $\vmax$ to capture the uncertainty of the neighboring agents in \eqref{eq:unc1}, a triggered event may occur faster than desired. Instead, the agents assume their neighbors are moving with a constant speed $\vpni$ when they are in motion. Let $\vpni$ be the maximum speed of the $i$th agent and its potential neighbors such that
\begin{align}\label{eq:vdg}
\vpni=\max\{\vi,\underset{j\in\pni}{\max(\vij)}\},
\end{align}
and let us rewrite \eqref{eq:unc1} in terms of $\vpni$ as following:
 
\begin{align}\label{eq:unc2}
\overline{\xij}=\OB \left(\pij,\left(\vpni (t-\tii)+
\vij (\tii -\tij)\right)\right).
\end{align}

\begin{remark}
The asymptotic convergence properties of the \algo~holds for any speed policy since the speed only effects the convergence time and nothing else. 
\end{remark}

Algorithm~\ref{tab:ecl1} summarizes the event-triggered broadcasting control algorithm with variable speed.

\vspace*{-0ex}
\begin{algorithm}[htb]
 {\footnotesize Initialization at time $t=0$, agent~$i\in \ID$ performs:
  \null\hfill\null \vspace*{-2ex}
  \begin{algorithmic}[1]
  \STATE sets $D =\pid\cup D_i$
  \STATE sets $\beta_i=1$
  \end{algorithmic}
}
  {\footnotesize At any time $t>0$, agent~$i\in \ID$ performs:
  \null\hfill\null \vspace*{-0ex}
  \begin{algorithmic}[1]
  \STATE updates $\pni$ using \eqref{eq:updateJ} 

  \IF{$\vi= 0$} 
  \STATE update $D =\pid\cup D_i$
  \STATE computes $\XX(D)$ as in \eqref{eq:unc}
  \ELSE
  \STATE computes $\overline{\XX}(D)$ as in \eqref{eq:unc2}
  \ENDIF
  \IF{$\vi = 0$ \textbf{ and }condition \eqref{Condition} is valid} 
  \STATE sets $\vi=\beta_i\vmax$
  \STATE broadcasts $\di$ using \eqref{eq:rad} $\ri(\JJ_i)$ distance away 
  \ELSIF{$\vi\neq 0 \textbf{ and }$condition \eqref{Condition} is invalid} 
  
  \STATE update $D =\pid\cup D_i$  
  \STATE computes $\XX(D)$ as in \eqref{eq:unc}
  \IF{$t-\tii<\Delta_{TB}$}
  \STATE sets $\beta_i=\beta_i/2$
  \ENDIF
  \IF{condition \eqref{Condition} is valid} 
  \STATE sets $\vi=\beta_i\vmax$
  \STATE broadcasts $\di$ using \eqref{eq:rad} $\ri(\JJ_i)$ distance away 
  \ELSE
  \STATE sets $\vi=0$ 
  \STATE broadcasts $\di$ using \eqref{eq:rad} $\ri(\JJ_i)$ distance away 
  \STATE waits for a time duration $\tau_d$
  \ENDIF
  \ENDIF

  \IF{$\big(t\neq t_i^i\big)$ \textbf{and} $\big($agent~$i$ has new potential neighbor \textbf{or} $\vi>\vii\big)$} 
  \STATE broadcasts $\di$ using \eqref{eq:rad} $\ri(\JJ_i)$ distance away 
  \ENDIF

  \STATE compute $u_i$ as in \eqref{eq:ui}
  \IF{$t-\tii>\Delta_{TB}$ \textbf{and} $\vi\neq 0$}
  \STATE sets $\beta_i=\min(2\beta_i,1)$
  \STATE sets $\vi=\beta_i\vmax$
  \ENDIF
  \end{algorithmic}
 }
 \caption{\hspace*{-.5ex}: \small \algo~ with variable speed}\label{tab:ecl1}
\end{algorithm}\vspace*{-0ex}

\section{Simulation}
In this section, we provide simulation results for the \algo. The simulations were developed using \texttt{MATLAB 2019a}. For all simulations, a time-step of $\step=1/60s$ was chosen {as if agents operating frequency is 60Hz} and all simulations were performed with eight agents~$N=8$ in a $40m\times 40m$ square environment. Agents were initialized with the locations, 
\[
\begin{split}
P=\big\{ (11.8, 36.3), (1.1, 6.0), (11.7,20.1),(15.3, 5.5), \\
 (11.6, 1.0), (7.5,9.1), (17.0, 15.3), (13.5, 6.3) \big\}.
\end{split}
\]
In addition, the initial value of the speed adjustment parameter $\beta_i$ and the maximum speed $\vmax$ were set to $\beta_i=1$ and $\vmax=0.1 m/s$ for all agents, respectively. The density function $\phi(q)$, { is provided} to all agents and was chosen to be $\phi(q)=e^{{-\TM{x-q_1}/100}}+e^{{-\TM{x-q_2}/100}}$ where $q_1=(20,30)$ and $q_2=(30,10)$. 

\begin{figure*}[htb]
    \centering
     \subfigure[]{\includegraphics[trim={1.8cm 1cm 1cm 1cm},clip,width=4.35cm]{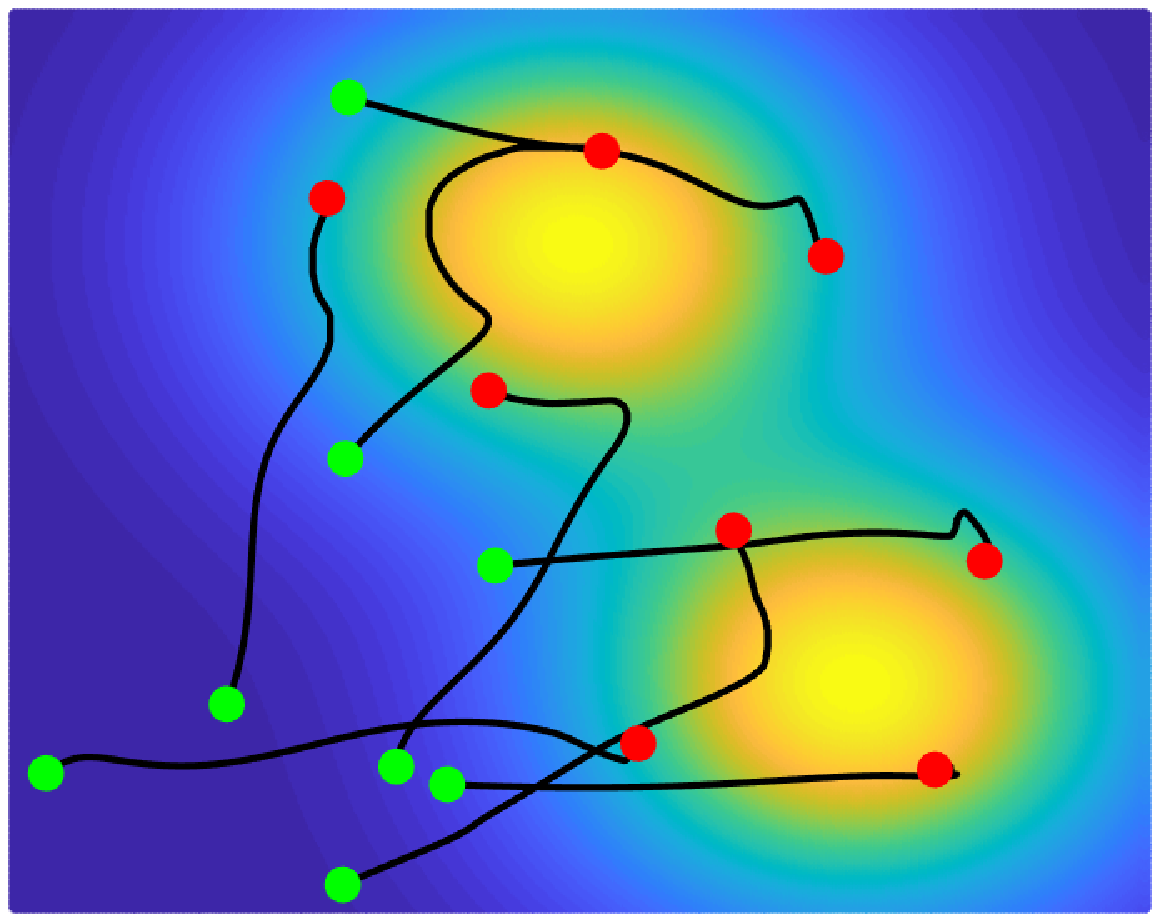}}
     \subfigure[]{\includegraphics[trim={1.8cm 1cm 1cm 1cm},clip,width=4.35cm]{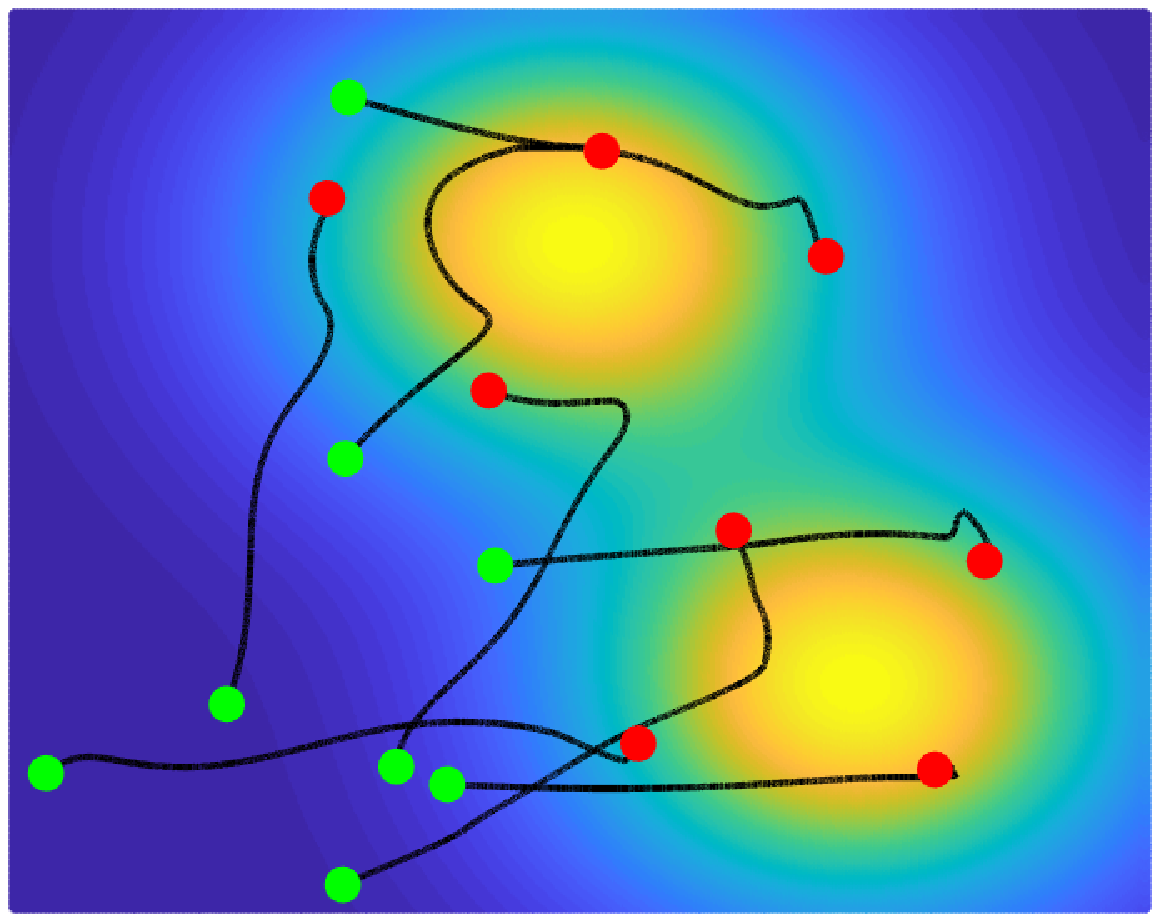}}
     \subfigure[]{\includegraphics[trim={1.8cm 1cm 1cm 1cm},clip,width=4.35cm]{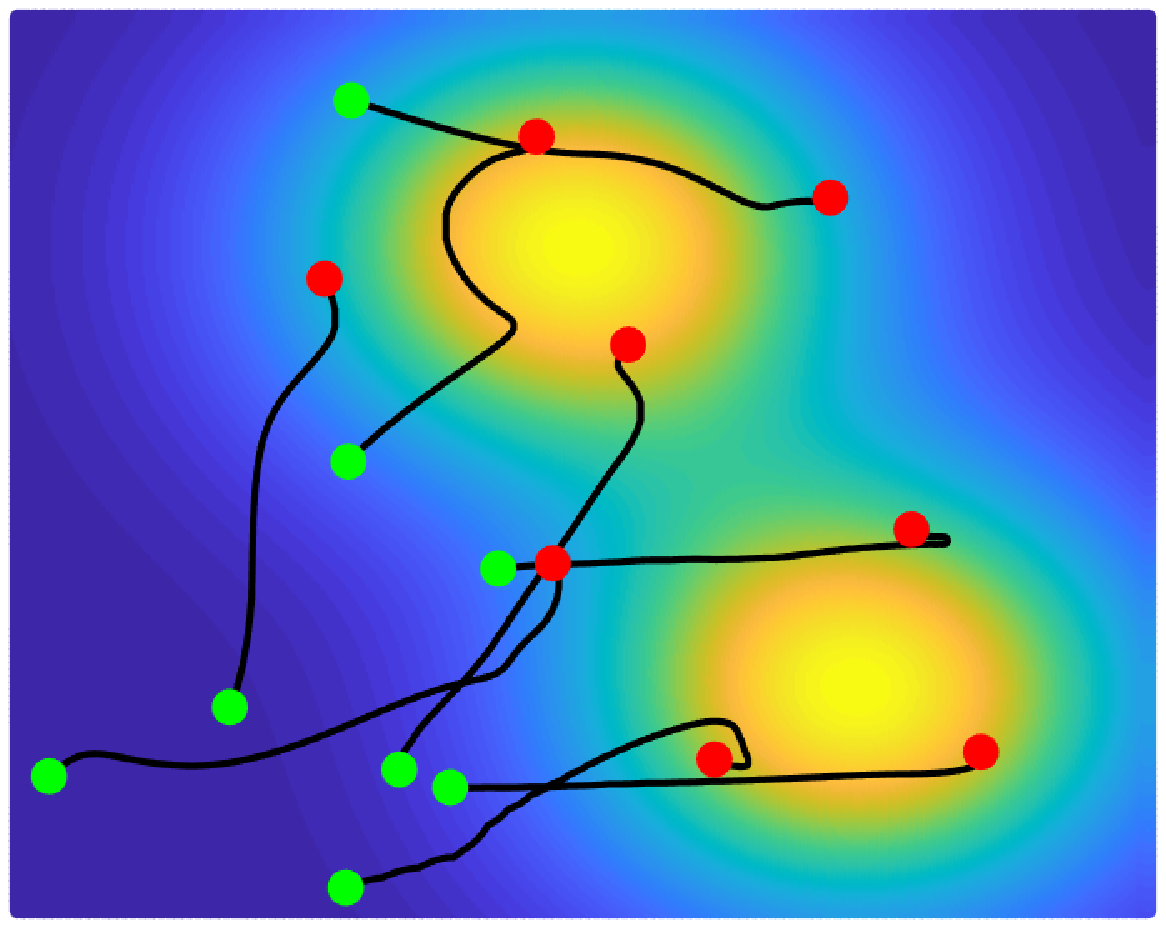}}     
     \subfigure[]{\includegraphics[trim={1.8cm 1cm 1cm 1cm},clip,width=4.35cm]{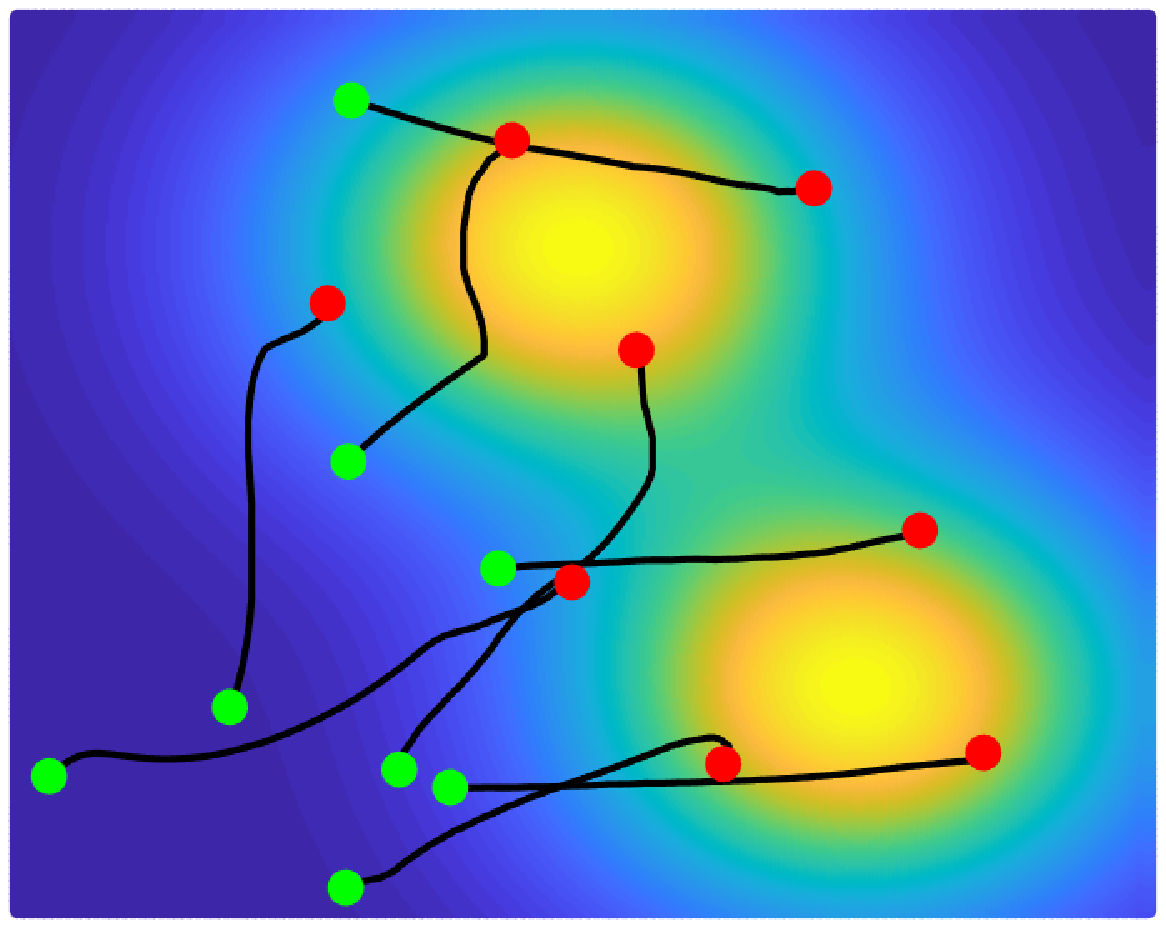}}
     \caption{Network trajectories of (a) periodic broadcasting algorithm,  (b) self-trigger algorithm in \cite{CN-JC:11-auto}, (c) \algo~constant $\vmax$ speed, and (d) \algo~with variable speed and $\Delta_{TB}=45/60s$. The green and red dots correspond to the initial and final agent positions, respectively}\label{fig:SI:D}
    \vspace*{-4ex}
\end{figure*}

\subsection{Simulation results}
The simulation results presented here demonstrate the effectiveness of the \algo~when compared to a periodic broadcasting algorithm where the agents broadcast at every time-step and compared to the self-triggered algorithm in \cite{CN-JC:11-auto}. According to Table 4, in \cite{CN-JC:11-auto} the self-triggered algorithm may require the agents to communicate more than once in the same time instance to ensure the sufficiency of the information. Instead for the self-triggered algorithm, we will assume the agents know their exact Voronoi neighbors at all times and only communicate once when a trigger is occurred.

\begin{figure}
\vskip3em
\psfrag{Periodic Broadcasting}{Periodic broadcasting algo}
\psfrag{Our algo with constant speed smax}{Our algo with constant speed}
\psfrag{Self}{Self-trigger algo}
\psfrag{TB=10/60}{Our algo with variable speed}
\psfrag{Objective Function Value}[Bc][0.9]{Objective function value $\HH$}
\psfrag{Time Step (s)}{Time (s)}
\psfrag{10}{$10$}
\psfrag{4}[Bcl][0.5]{$4$}

\psfrag{Y1}{$1$}
\psfrag{Y2}{$2$}
\psfrag{Y3}{$3$}
\psfrag{Y4}{$4$}
\psfrag{Y5}{$5$}
\psfrag{Y6}{$6$}
\psfrag{Y7}{$7$}
\psfrag{Y8}{$8$}
\psfrag{Y9}{$9$}
\psfrag{Y10}{$10$}
\psfrag{X0}{$0$}
\psfrag{X100}{$100$}
\psfrag{X200}{$200$}
\psfrag{X300}{$300$}
\psfrag{X400}{$400$}
\psfrag{X500}{$500$}
\psfrag{X600}{$600$}
\psfrag{X50}{$50$}
\psfrag{X150}{$150$}
\psfrag{X250}{$250$}
\psfrag{X350}{$350$}
\psfrag{X450}{$450$}
\psfrag{X550}{$550$}
\includegraphics[width=8.3cm]{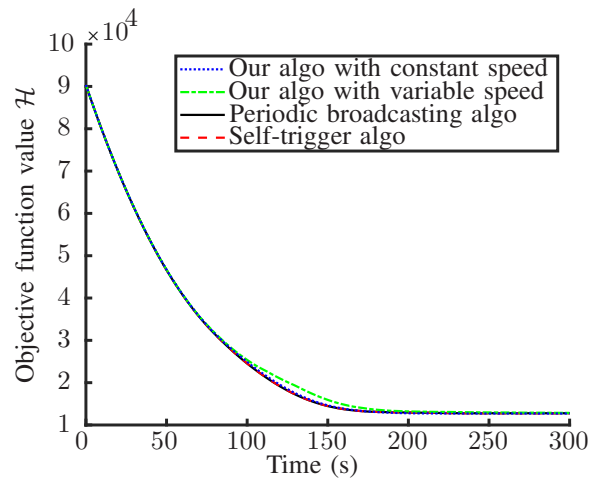}
\centering
\caption{A comparison of the objective function value}\label{fig:SI:H}
\end{figure}

Let us start with Fig.~\ref{fig:SI:H}. This figure shows the algorithms comparison for the convergence of the objective function. {It can be seen in Fig. 4 that all algorithms including our algorithm with constant speed $\vmax$, with variable speed and $\Delta_{TB}=45/60s$, the self-trigger algorithm, and the periodic broadcasting algorithm reach similar objective function values. This implies that the reduced communication by our algorithm does not affect the convergence to the set of centroidal Voronoi configurations. In other words, our \algo~ preforms as good as the periodic broadcasting algorithm with respect to $\HH$. } 

\begin{figure}
\vskip3em
\psfrag{Periodic Broadcasting}{Periodic broadcasting algo}
\psfrag{Our algo with constant speed smax}{Our algo with constant speed}
\psfrag{Self}{Self-trigger algo}
\psfrag{TB=10/60}{Our algo with variable speed}
\psfrag{Amount of Communication}[Bc][0.9]{Amount of communication}

\psfrag{Time}{Time (s)}
\psfrag{10}{$10$}
\psfrag{3}[Bcl][0.5]{$3$}

\psfrag{Y0}{$0$}
\psfrag{Y20}{$20$}
\psfrag{Y40}{$40$}
\psfrag{Y60}{$60$}
\psfrag{Y80}{$80$}
\psfrag{Y100}{$100$}
\psfrag{Y120}{$120$}
\psfrag{Y140}{$140$}
\psfrag{Y160}{$160$}
\psfrag{Y180}{$180$}
\psfrag{X0}{$0$}
\psfrag{X100}{$100$}
\psfrag{X200}{$200$}
\psfrag{X300}{$300$}
\psfrag{X400}{$400$}
\psfrag{X500}{$500$}
\psfrag{X600}{$600$}
\psfrag{X50}{$50$}
\psfrag{X150}{$150$}
\psfrag{X250}{$250$}
\psfrag{X350}{$350$}
\psfrag{X450}{$450$}
\psfrag{X550}{$550$}
\includegraphics[width=8.3cm]{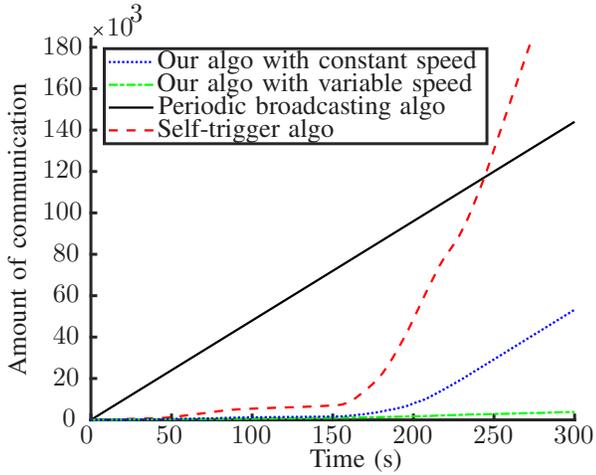}
\centering
\caption{A comparison of the amount of communications between the agents}\label{fig:SI:C}
\end{figure}

Additionally, Fig.~\ref{fig:SI:C} shows the algorithms comparison for the amount of communication between agents. {It is clear that our algorithm with constant speed of $\vmax$ and our algorithm with variable speed with $\Delta_{TB}=45/60s$ both significantly reduced the amount of communication between the agents. It is clear that the self-triggered algorithm performs worse than the periodic broadcasting algorithm due to the fact that the self-triggered algorithm requires agents to both request and respond. The amount of communication is greater than the periodic broadcasting algorithm at the $300s$ mark. We would like to note that a request for information is counted as a single communication and each response is counted as a single communication as well. For example, if an agent sent a request and $5$ agents respond, this becomes a total of $6$ communicated messages. Now the advantage of relaxing the assumption of instantaneous communication and synchronous actions by agents is shown clearly in the comparison between the self-trigger algorithm and our \algo.}
To quantify our results, it is noted that for our algorithm with constant speed of $\vmax$, the amount of communication between agents is reduced by $63.0\%$ and $78.2\%$ when compared to the periodic broadcasting and self-triggered algorithms, respectively. For our algorithm with variable speed and with $\Delta_{TB}=45/60s$, the amount of communication between agents is reduced by $97.3\%$ and $98.4\%$ when compared to the periodic broadcasting and self-triggered algorithms, respectively.

\begin{figure}
\vskip3em
\psfrag{Amount of Communication}[Bc][0.9]{Amount of communication per time step}
\psfrag{Time}{Time (s)}
\psfrag{X0}{$0$}
\psfrag{X100}{$100$}
\psfrag{X200}{$200$}
\psfrag{X300}{$300$}
\psfrag{X400}{$400$}
\psfrag{X500}{$500$}
\psfrag{X600}{$600$}
\psfrag{X50}{$50$}
\psfrag{X150}{$150$}
\psfrag{X250}{$250$}
\psfrag{X350}{$350$}
\psfrag{X450}{$450$}
\psfrag{X550}{$550$}
\psfrag{Y1}[Bc][0.9]{$1$}
\psfrag{Y2}[Bc][0.9]{$2$}
\psfrag{Y3}[Bc][0.9]{$3$}
\psfrag{Y4}[Bc][0.9]{$4$}
\psfrag{Y5}[Bc][0.9]{$5$}
\psfrag{Y6}[Bc][0.9]{$6$}
\psfrag{Y7}[Bc][0.9]{$7$}
\psfrag{Y8}[Bc][0.9]{$8$}
\psfrag{Y9}[Bc][0.9]{$9$}
\psfrag{Y0}[Bc][0.9]{$0$}
\includegraphics[width=8.3cm]{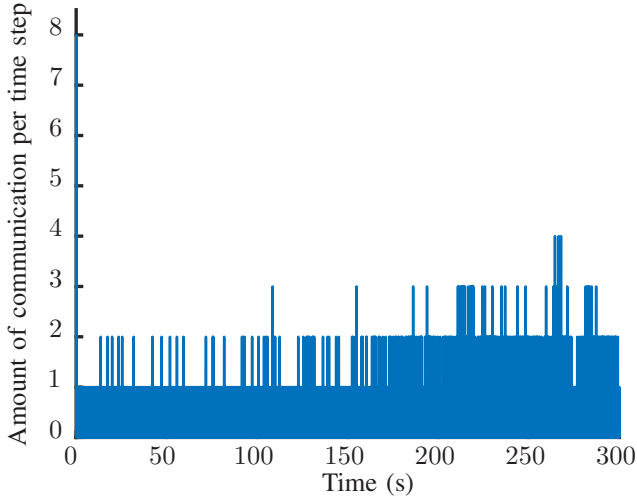}
\centering
\caption{The amount of communications per time step under the \algo~with variable speed and $\Delta_{TB}=45/60$}\label{fig:SI:D1}
\end{figure}

Last but not least, our algorithm with constant $\vmax$ speed (around the $230$s mark) will require periodic communication since condition \eqref{Condition} becomes invalid after a time-step for all agents. In addition, our algorithm with variable speed can be seen to address the issue of eventual constant communication as shown in Fig.~\ref{fig:SI:D1} and was able to further reduce the amount of communication as shown in Fig.~\ref{fig:SI:C}. It is notable that in order to see the effectiveness of the \algo~with variable speed, the target time duration $\Delta_{TB}$ must be greater than or equal to twice the time-step, $\Delta_{TB}\geq 2\step$. In other words, when $\Delta_{TB}\geq 2\step$, the agents find the appropriate speeds such that they can move at least for $2\step$ without broadcasting. Also as $\Delta_{TB}$ increases, the amount of communication is further reduced with a small delay in convergence speed as is seen in Fig.~\ref{fig:SI:C} and Fig.~\ref{fig:SI:H} for $\Delta_{TB}=45/60s$.

In summary, these figures illustrate how the \algo~is able to achieve similar convergence performance to both the periodic broadcasting and self-triggered algorithm while also requiring much less communication between agents. As the figures show, there is a trade-off between the convergence speed and the amount of communication that occurs between agents. With slightly slower convergence speed, our algorithm significantly reduced the amount of communication compared to the periodic broadcasting and the self-triggered cases. The reader may note that given the communication range assignment for the sefl-triggered algorithm in \cite{CN-JC:11-auto}, the algorithm will require much more communication than showing in Fig.~\ref{fig:SI:C}.

\section{Conclusion and Future Work}
{This paper proposed a distributed \algo~implemented with a more practical and fully asynchronous broadcasting communication model. The broadcasting communication model relaxes the assumptions of instantaneous communication and synchronous actions by agents. In other words, this means that the communication model does not allow agents to request information from other agents nor acknowledge the reception of messages. Instead, agents must strictly broadcast their information when they decide that it is appropriate.}
In addition, the agents are capable of determining, prior to transmission, the sufficient broadcasting range to share their information with their potential neighbors i.e. all Voronoi neighbors.
Through analysis, the proposed algorithm was shown to provide guaranteed asymptotic convergence. Also, the algorithm was shown to significantly reduce the amount of communication between agents when compared with both a periodic broadcasting strategy and a self-triggered request-response strategy. Future work will focus on a modification of the broadcast deployment problem where the possibility of packet drops may occur. We aim to use the possible and potential neighbors and the communication range assignment method presented here to allow agents to not only determine when to communicate, but also how to select a channel to communicate over in order to minimize packet loss.

\appendices 
\newcommand{\dgvji}{dgV_{i,j}}
\section{Proof of Lemma~\ref{lm:ndg}}\label{app:A}
\begin{IEEEproof}  
For convenience, let $\dgvji = \left(\dgvi(\XX) \cap \dgvj(\XX)\right)$ and let $\ID_{dg}$ be the set of agents such that $\dgvji \cap \dgvk \neq \emptyset$. Now, if $\dgvji\neq \emptyset$, we know 
\begin{align*}
\overline{\dgvji}=\dgvji  \setminus \cup_{k \in \ID_{dg}}\gvk(\XX)\neq \emptyset
\end{align*}
since $\gvi(\XX)\subset \dgvi(\XX)\;\forall i\in\ID_{dg}$.
The point in $\overline{\dgvji}$ are not guaranteed to be closer to any agent in $\ID_{dg}$ than the others, and they are guaranteed to be closer these agents than agents $k\in\ID\setminus \ID_{dg}$ since it is outside all the others' dual-guaranteed cells. 
Therefore, $\exists P\subset \XX$ such that agent~$i$ and $j$ share at least one point in $\overline{\dgvji}$ such that $\voi\cap \voj=\{q\}\;|\;q\in \overline{\dgvji}$.
\end{IEEEproof}

\section{Proof of Corollary~\ref{coro}}\label{app:B}
\begin{IEEEproof}  
Let us consider the worst case for $\dgvi(\XX) \cap \dgvj(\XX) \neq \emptyset$ to be true that is
\begin{align*}
\dgvi(\XX) \cap \dgvj(\XX)=\{q'\}\;|\; q'=\underset{q\in \dgvi(\XX)}{\max}(\TM{q-\ppi}).
\end{align*}
Given the definition of the dual-guaranteed Voronoi cell, this implies that $\min\TM{q'-\xli} = \max\TM{q'-\xlj}\; \forall \xli\in \xi^i, \;\xlj\in \xj^i$. 
Since agent~$i$ knows its exact location, we can rewrite the previous equation as 
\begin{align*}
\TM{q'-\ppi} = \underset{\xlj\in \xij}{\max}\TM{q'-\xlj}.
\end{align*} 
By rearranging the previous equation, we get
\begin{align*}
\underset{\xlj\in \xij}{\max}\TM{\ppi-\xlj} \leq 2\TM{q'-\ppi}.
\end{align*}
Thus, any agent $j$ s.t. $\pj \not\in \overline{B}(\ppi,2\TM{q'-\ppi})$ is guaranteed not to be a dual-guaranteed neighbor since the dual-guaranteed Voronoi cells will not intersect and is guaranteed not to be a Voronoi neighbor by Lemma~\ref{lm:ndg}.
\end{IEEEproof}

\section{Proof of Proposition~\ref{pro:newBound}}\label{app:E}
{
\begin{IEEEproof}  
To prove the claim we must show that both
\begin{align}\label{eq:show1}
\TM{\cdgi-C_{V_i}} \leq \bnd 
\end{align}
and
\begin{align}\label{eq:show2}
\TM{\cddgi-C_{V_i}} \leq \bnd
\end{align}
hold. By~\cite[Proposition~5.2]{CN-JC:11-auto},
we know that for any sets~$L \subset V \subset U$, 
\begin{align*}
\TM{C_V-C_L}\leq 2cr_{U} \left(1- \frac{M_{L}}{M_{U}} \right).
\end{align*}
Since~$\gvii(\dgi) \subset \voi \subset \dgvii(\dgi)$,
the first condition~\eqref{eq:show1} follows immediately with~$L = \gvii, V = V_i,$ and~$U = \dgvii$. To show~\eqref{eq:show2}, let~$L = V_i,$ and~$V=U=\dgvii$, then
\begin{align*}
\TM{\cddgi-C_{V_i}} &\leq 2cr_{\dgvii} \left(1- \frac{M_{V_i}}{M_{\dgvii}} \right) \\
&\leq 2cr_{\dgvii} \left(1- \frac{M_{\gvii}}{M_{\dgvii}} \right) = \bnd, 
\end{align*}
which concludes the proof. 

\end{IEEEproof}  
}
\section{Proof of Lemma~\ref{lm:J}}\label{app:C}
\begin{IEEEproof}  
We start by proving $j\in\pni\Leftrightarrow i\in\pnj$ at all times. 
For $j\in\pni\Leftrightarrow i\in\pnj$ to holds, agent $i$ and $j$ to must have the same last broadcasted information from their common potential neighbors and their-self. The reason begin that any agent that is not a common neighbor cannot be closer to any point in $\dgvi(\XX)\cap\dgvj(\XX)$ than agent $i$, $j$ and their common neighbors. 
Under the broadcasting range assignment, each common potential neighbor's information is guaranteed to reach agent $i$ and $j$. In addition, by \eqref{eq:updateJ}, if $j\in\pni$, agent $i$ uses to the same information agent $j$ has about it-self $\xii$ and agent $j$ will do the same as well. In case $j\not\in\pni$, agent $i$ uses its prefect information since it will broadcast to agent $j$ by the decision control law in Section \ref{se:DC} if they becomes new neighbor, and when agent $j$ receives agent $i$'s information, agent $j$ will have agent $i$'s prefect information. 
Therefore, agent $i$ and $j$ will always have the same information required to determine if they are or they are not potential neighbor that guarantees $j\in\pni\Leftrightarrow i\in\pnj$ at all times.

Now, we want to show $\nei\subset \pni$ is guaranteed at event-times. When agent $i$ broadcast $\ri(\pni)$ distance away using the \eqref{eq:rad}, by Corollary \ref{coro}, all its potential neighbors will receive its broadcasted. Since $j\in\pni\Leftrightarrow i\in\pnj$ and by the decision control law in Section \ref{se:DC}, any agent $k\not\in\pni$ will broadcast its information as soon as it gets agent $i$ information if they become potential neighbors. Thus, when agent $i$ broadcast, it will receive the new potential neighbors' information immediately. In case the $i$th agent did not receive any information when it broadcast, it implies that the agent does not have any new neighbor. Therefore, $\nei\subset \pni$ is guaranteed at event-times.
\end{IEEEproof}

\section{Proof of Proposition~\ref{pro:algo}}\label{app:D}
\begin{IEEEproof}  
In this proof, we want to guarantee that 
\begin{align}\label{APD:1}
\gvii(\pni)\subset \voi\subset \dgvii(\pni)
\end{align} 
at all times under the the \algo. Let us start by saying that as the potential neighbors' uncertainties increase, the guaranteed Voronoi cell shrinks and the dual-guaranteed Voronoi cell expands, and when the agent computes the uncertainties using \eqref{eq:unc1}, the guaranteed and dual-guaranteed Voronoi cells change faster than when the uncertainties computed by \eqref{eq:unc}.
By~\cite[Lemma~4.1~and~4.2]{CN-JC:11-auto}, 
the $i$th agent' cells given the agent's and $\pid$ information satisfies $\gvii(\overline{\XX})\subset \gvii(\XX)\subset \voi\subset \dgvii(\XX)\subset \dgvii(\overline{\XX})$, where $\overline{\XX}=\{\xij\}_{j\in\pni}$ computed using \eqref{eq:unc1} and $\XX=\{\xij\}_{j\in\pni}$ computed using \eqref{eq:unc}. 
Let us start by proving \eqref{APD:1} is guaranteed at every event-time. By Lemma~\ref{lm:J}, at every event-time, $\nei\subset\pni$ is guaranteed, and as a result \eqref{APD:1} is guaranteed as well since the $i$th agent has all Voronoi neighbors' information. 

Now, we will prove that \eqref{APD:1} is guaranteed between event-times. Ideally, the $i$th agent can move and compute the uncertainties using \eqref{eq:unc} until condition \eqref{Condition} is invalid. However, this requires $\nei\subset \pni$ to be true at all times. This is challenging to ensure because agent $i$ will not know about a new potential neighbor until it broadcast. Instead, we let the agents compute the uncertainties using \eqref{eq:unc1} when they are in motion.
\cite[Lemma~4.1]{CN-JC:11-auto} state that if $\nei\subset \pni$ it satisfied and the agent expands the uncertainties using \eqref{eq:unc1}, $\gvii(\pni)\subset \voi$ is guaranteed without using any additional information.
In addition, \cite[Lemma~4.2~and~4.3]{CN-JC:11-auto} state that by expanding $\{\xij\}_{j\in\pni}$ using \eqref{eq:unc1}, the dual-guaranteed Voronoi cell cannot be bigger given any agent $k\in\ID$ perfect information. In fact \cite[Lemma~4.2~and~4.3]{CN-JC:11-auto} guarantee $\voi\subset \dgvii(\pni)$ at all times even if $\nei\not\subset\pni$. Since $\nei\subset\pni$ is guaranteed at event-time, by Lemma~\ref{lm:J}, under the \algo, \eqref{APD:1} is guaranteed while the agents are in motion.  

Furthermore, when the $i$th agent is waiting for new/updated information, the agent will not affect any agent $k\in\ID$ because it is not moving. Also,
since \eqref{APD:1} is guaranteed while the other agents are moving, the moving agents will not affect agent $i$. Therefore, agent $i$ cannot affect or be affected by any other agent. Thus, the agent expands the uncertainties as necessary using \eqref{eq:unc}, and the \algo~guaranteed \eqref{APD:1} while the agents are waiting.  

Since \eqref{APD:1} is guaranteed at event-times and between event-times, \eqref{APD:1} is guaranteed at all-times.
\end{IEEEproof}

\bibliographystyle{IEEEtran}
\bibliography{Paper}

\begin{IEEEbiography}[{\includegraphics[width=1in,height=1.25in,clip,keepaspectratio]{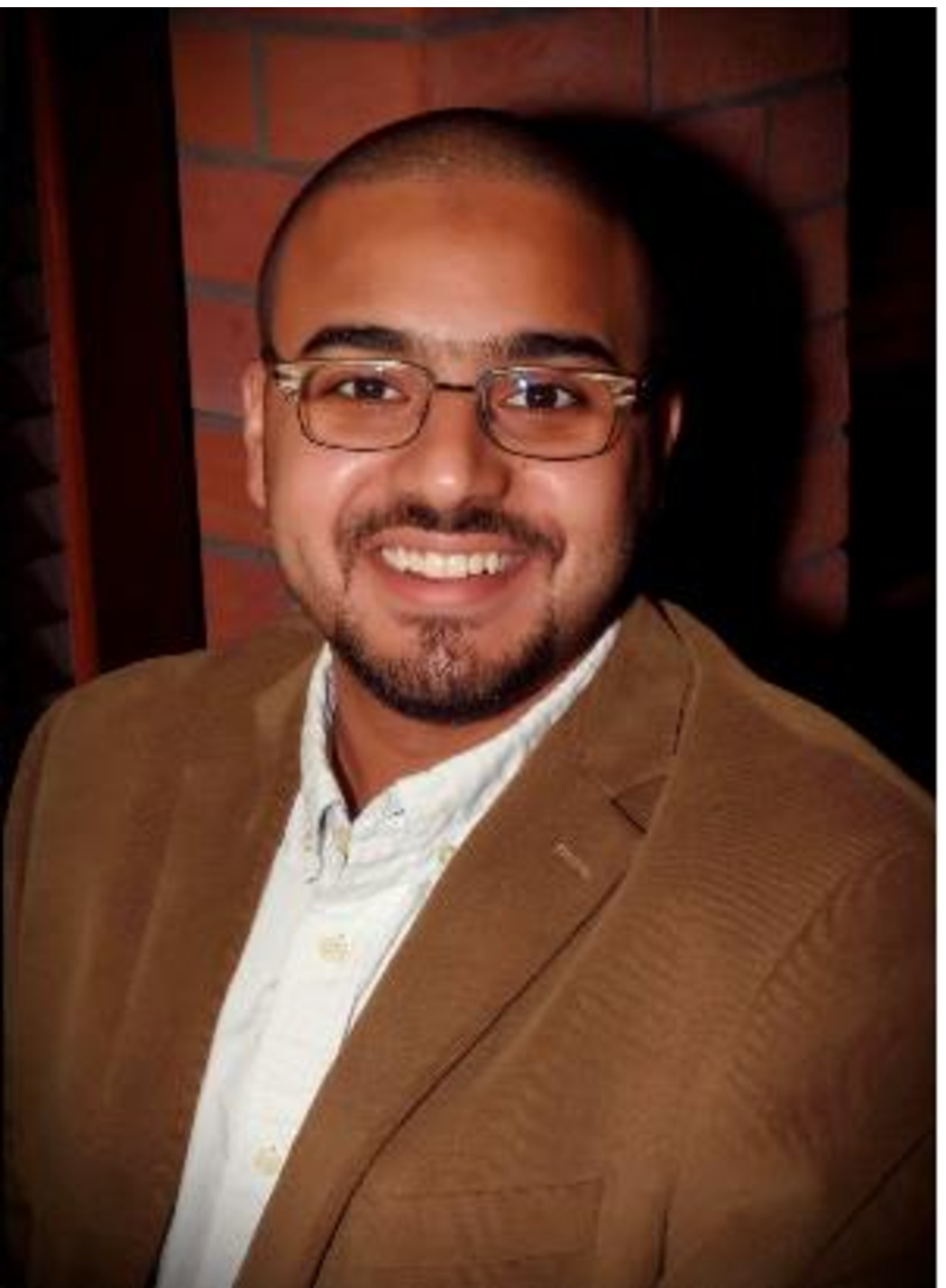}}]{Mohanad~Ajina} is PhD candidate in Electrical and Computer Engineering department at George Mason University. He received his undergrad from Southern Illinois University Carbondale in 2014. In th1e same year, he joined George Mason University to do his Master in Electrical Engineering and he graduated in 2016. His current research interests include distributed coordination algorithms, robotics, event- and self-triggered control, advanced data analytic, statistical modeling and machine learning.
\end{IEEEbiography}

\begin{IEEEbiography}[{\includegraphics[width=1in,height=1.25in,clip,keepaspectratio]{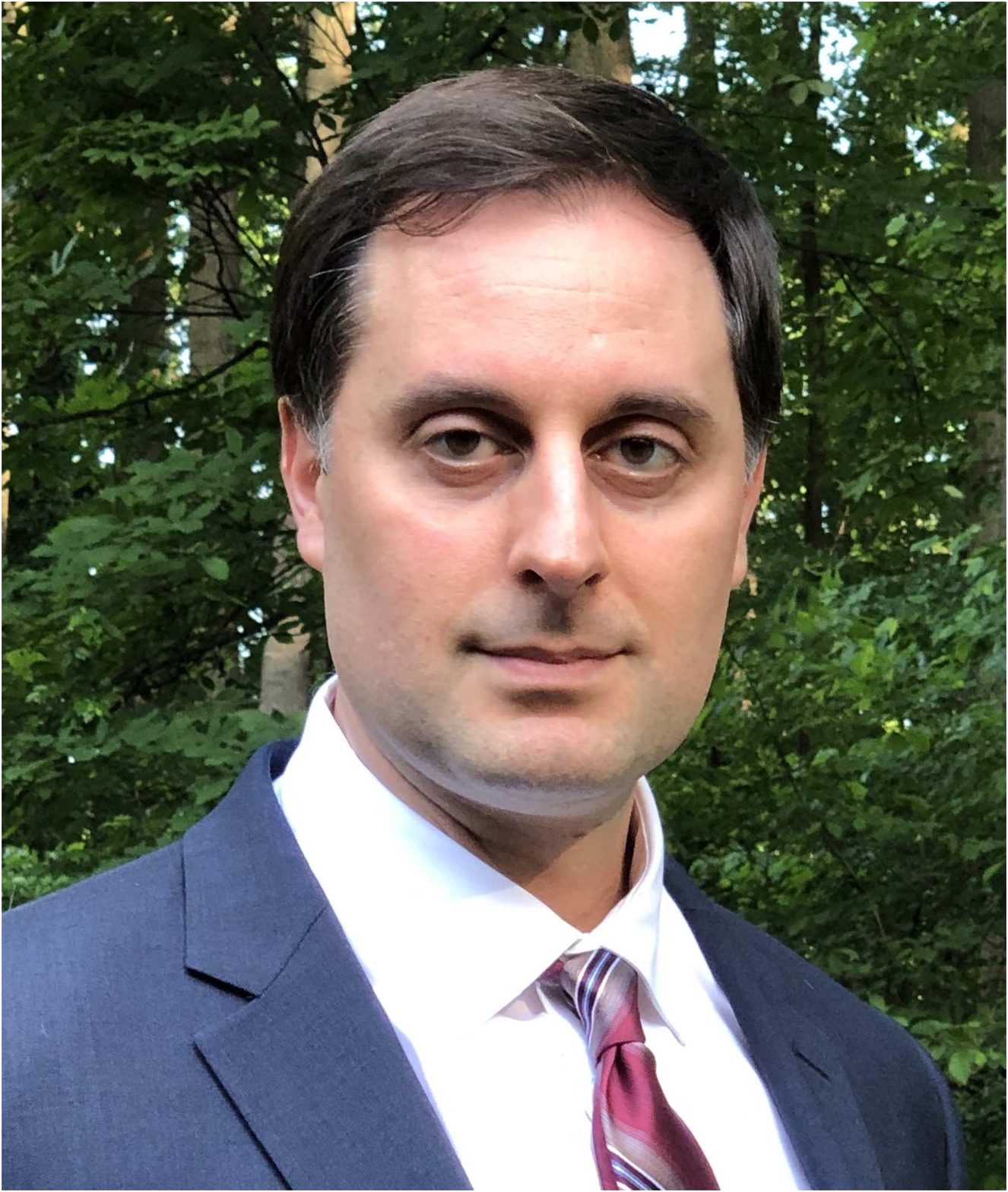}}]{Daniel~Tabatabai} is a PhD student in the Electrical and Computer Engineering department at George Mason University and a Senior Robotics Engineer for d3i Systems. He received both the B.S. and M.S. degrees in Electrical Engineering from George Mason University in 2004 and 2014, respectively. His research interests include distributed coordination algorithms, event- and self-triggered control, autonomous robotic locomotion and control, optimization and reinforcement learning in control, and physiological modeling and control.
\end{IEEEbiography}

\begin{IEEEbiography}[{\includegraphics[width=1in,height=1.25in,clip,keepaspectratio]{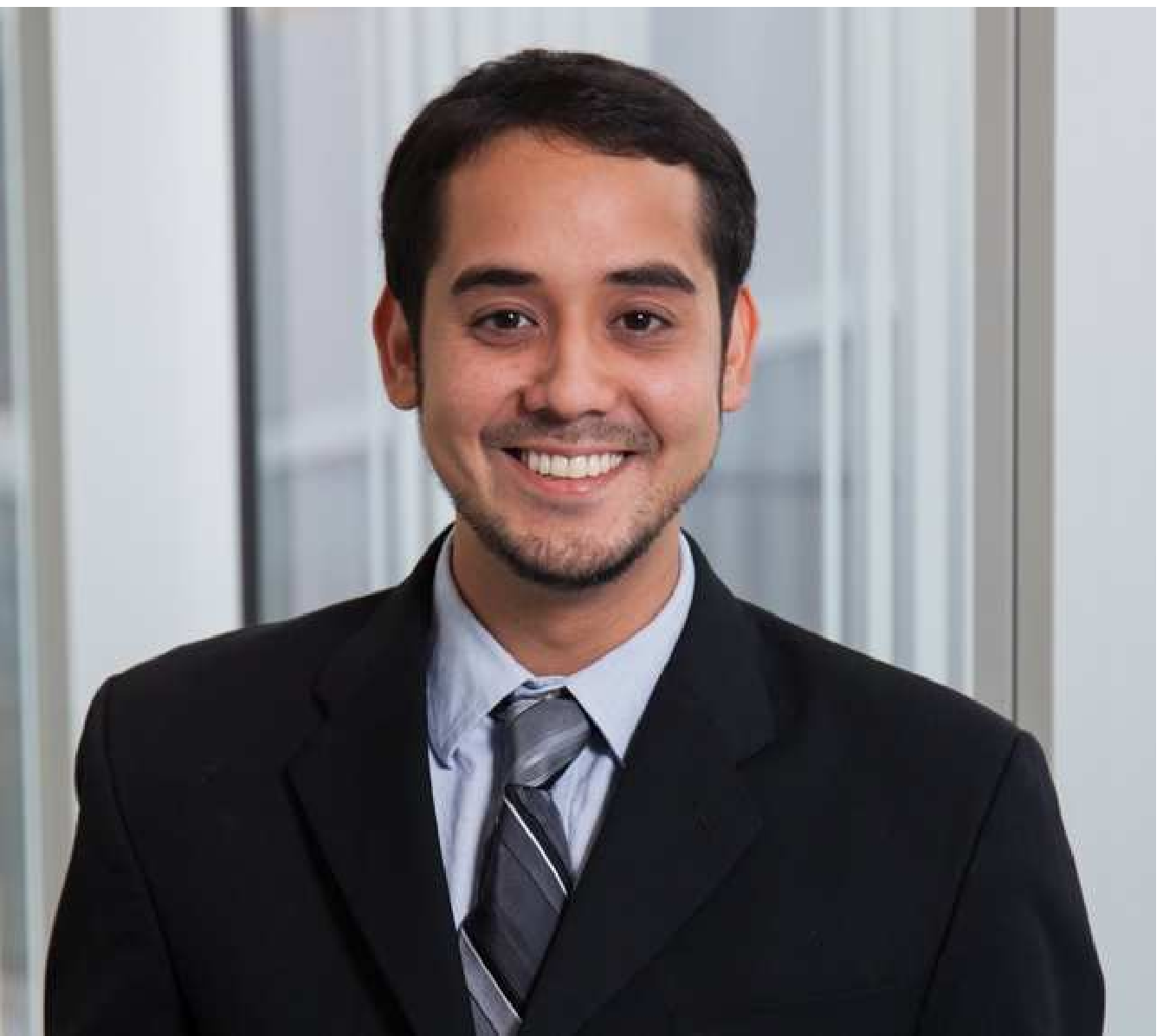}}]{Cameron~Nowzari} received the Ph.D. in Mechanical Engineering from the University of California, San Diego in September 2013. 
He then held a postdoctoral position with the Electrical and Systems Engineering Department at the University of Pennsylvania until 2016. He is currently an Assistant Professor with the Electrical and Computer Engineering Department at George Mason University, in Fairfax, Virginia. He has received
several awards including the American Automatic Control Council's O. Hugo Schuck Best Paper
Award and the IEEE Control Systems Magazine Outstanding Paper Award. His current research interests include dynamical systems and control, distributed coordination algorithms, robotics, event- and self-triggered control, Markov processes, network science, spreading processes on networks and the Internet of Things.
\end{IEEEbiography}

\end{document}